\documentclass[12pt]{amsart}
\usepackage{mathrsfs, amssymb,  euscript}%, diagrams}
\usepackage[matrix,arrow,curve]{xy}
\makeatletter \@addtoreset{equation}{section} \makeatother

\renewcommand\labelenumi{(\roman{enumi})}
\newcommand{\xref}[1]{{\rm \ref{#1}}}
\newcommand{\comment}[1]{}
\newcommand{\comp}{\mathbin{\scriptstyle{\circ}}}

\newcommand{\bir}{\approx}
\newcommand{\Aut}{\operatorname{Aut}}
\newcommand{\Sing}{\operatorname{Sing}}

\newcommand{\St}{\operatorname{St}}
\newcommand{\Pic}{\operatorname{Pic}}
\newcommand{\Br}{\operatorname{Br}}
\newcommand{\Hom}{\operatorname{Hom}}
\newcommand{\Mori}{\overline{NE}}

\newcommand{\point}[2]{P\bigl[^{#1}_{#2}\bigr]}
\newcommand{\plane}[2]{\Pi\bigl[^{#1}_{#2}\bigr]}

\newcommand{\chit}{\chi_{\operatorname{top}}}

\newcommand{\Char}{\operatorname{char}}
\newcommand{\Diag}{\operatorname{Diag}}

\newcommand{\CC}{{\mathbb C}}

\newcommand{\ZZ}{{\mathbb Z}}
\newcommand{\QQ}{{\mathbb Q}}
\newcommand{\PP}{{\mathbb P}}
\newcommand{\FF}{{\mathbb F}}

\newcommand{\II}{{\mathbb I}}
\newcommand{\OO}{{\mathbb O}}
\newcommand{\TT}{{\mathbb T}}

\newcommand{\AAA}{{\mathfrak A}}
\newcommand{\SSS}{{\mathfrak S}}
\newcommand{\DDD}{{\mathfrak D}}

\newcommand{\OOO}{{\mathscr{O}}}

\newcommand{\EEE}{{\mathscr{O}}}

\newcommand{\muu}{{\boldsymbol\mu}}
\newtheorem{theorem}[equation]{Theorem}
\newtheorem{proposition}[equation]{Proposition}
\newtheorem{lemma}[equation]{Lemma}
\newtheorem{corollary}[equation]{Corollary}

\theoremstyle{definition}
\newtheorem{definition}[equation]{Definition}
\newtheorem{example}[equation]{Example}
\newtheorem{problem}[equation]{Problem}
\newtheorem{emptytheorem}[equation]{}

\newtheorem{remark}[equation]{Remark}

\newcounter{ssylki}[equation]
\renewcommand{\thessylki}{(\arabic{section}.%
\arabic{equation}.\arabic{ssylki})}

\title{Fields of invariants of finite linear groups}
\author{Yuri G. Prokhorov}
\address{
Department of Algebra, Faculty of Mathematics (and Mechanics),
Moscow State University, Moscow 117234, Russia}
\email{prokhoro@mech.math.msu.su}

\thanks{The work was partially supported by grants
CRDF-RUM, No. 1-2692-MO-05 and
RFBR, No. 05-01-00353-a, 06-01-72017.
}

\begin{document}
\maketitle
\section{Introduction}
Let $G$ be a finite group of order $n$ and let $\Bbbk$ be a field.
Consider a rational (i.e., pure transcendental) extension
$K/\Bbbk$ of transcendence degree $n$. We may assume that
$K=\Bbbk(\{x_g\})$, where $g$ runs through all the elements of
$G$. The group $G$ naturally acts on $K$ via $h(x_g)=x_{hg}$. E.
Noether \cite{Noether_1913}
asked whether the field of invariants $K^G$ is rational over
$\Bbbk$ or not. On the language of algebraic geometry, this is 
a question about the rationality of the
quotient variety $\mathbb A^n/G$.

The most complete answer on this question is known for Abelian
groups. Thus, if $G$ is Abelian of exponent $e$, $\Char \Bbbk$
does not divide $e$ and $\Bbbk$ contains a primitive $e$th roots
of unity then $\mathbb A^n/G$ is rational \cite{Fischer-E-1915-a}. 
On the other hand, over an arbitrary
field $\Bbbk$ the rationality of $\mathbb A^n/G$ is related to
some number theoretic questions. In this case $\mathbb A^n/G$ can
be non-rational \cite{Swan-1969} (see also
\cite{Voskresenski-1973-re}, \cite{Endo-Miyata-1973},
\cite{Lenstra-1974}).

Noether's question can be generalized as follows.
\begin{problem}
\label{problem-main}
Let $G$ be a finite group, let $V$ be a finite-dimensional vector
space over an algebraically closed field $\Bbbk$, and let
$\rho\colon G\to GL(V)$ be a representation. Whether the quotient
variety $V/G$ is rational over $\Bbbk$?
\end{problem}

The following fact is well known to specialists.
\begin{proposition}
\label{prop-first}
$V/G$ is birationally equivalent to $\PP(V)/G\times\PP^1$.
\end{proposition}
Indeed, by blowing up the origin we get a generically 
$\CC^1$-fibration $f\colon \widetilde {V/G}\to \PP(V)/G$
which is locally trivial in the \'etale topology.
By the geometric form of 
of Hilbert's theorem 90 the fibration
$f$ also locally trivial
in the Zariski topology.
We refer to  \cite{Miyata-K-1971}, \cite{Kervaire-Vust-1989},
\cite{colliotthelene-Sansuc-2005} 
for detailed algebraic proofs.

By Proposition \ref{prop-first} the rationality of $\PP(V)/G$ implies the rationality of
$V/G$. The inverse implication is not known (cf.
\cite{Beauville-Colliot-Thelene-Sansuc-1985}). The affirmative
answer to the L\"uroth problem gives us the rationality of $\PP(V)/G$
(and therefore $V/G$) for $\dim V\le 3$ over any algebraically
closed field, cf. \cite[Ch. 17]{Burnside-1911}. Thus $\dim V=4$ is the
first nontrivial case.

Most of this survey is devoted to reviewing the current status of
this problem in the special case where $V$ is of dimension four.
We refer to surveys \cite{Formanek-1984}, \cite{Dolgachev-1987},
\cite{Kervaire-Vust-1989}, \cite{Hajja-2000-surv},
\cite{colliotthelene-Sansuc-2005}  for other
aspects of \xref{problem-main}. 
 
If the representation $G \hookrightarrow GL(V)$ is not
irreducible, then the decomposition $V=V_1\oplus V_2$ gives us a
$G$-equivariant $\PP^1$-bundle structure $\widetilde {\PP(V)}\to
\PP(V_1)\times \PP(V_2)$, where $\widetilde {\PP(V)}$ is the
blowup of $\PP(V)$ along $\PP(V_1)\cup \PP(V_2)$. Obviously, this
$\PP^1$-bundle has a $G$-invariant sections, exceptional divisors
over $\PP(V_1)$ and $\PP(V_2)$. Therefore, $\PP(V)/G\bir
(\PP(V_1)\times \PP(V_2))/G\times \PP^1$. If $\dim V=4$, this
implies the rationality of $\PP(V)/G$. Thus, in the rationality
question of $\PP^3/G$, we may assume that $G\subset GL(V)$ is
irreducible.

\begin{definition}
\label{def-imp}
Let $G\subset GL(V)$ be a finite irreducible subgroup. We say that
$G$ is \emph{imprimitive of type $(m^k)$} if $V$ is a nontrivial
direct sum $V=V_1\oplus \cdots\oplus V_k$ with $\dim V_i=m$ such
that $gV_i=V_j$ for all $g\in G$. If such a direct splitting does
not exists, $G$ is called \emph{primitive}.
\end{definition}

As a consequence of Jordan's theorem \cite{Jordan-1870} one can
see that modulo scalar matrices in any dimension $n$ there is only
a finite number of primitive finite subgroups $G\subset
GL(n,\CC)$. Lower-dimensional primitive groups have been completely classified, 
see references in \cite[\S 8.5]{Feit-1971}.

\begin{emptytheorem}
\label{subsect-def-groups}
Let $G\subset GL(4,\CC)$ be a finite subgroup. 
By \xref{prop-first} we may regard $G$ modulo scalar matrices. 
Following \cite{Blichfeldt} we reproduce the
list of all primitive subgroups $G\subset GL(4,\CC)$ modulo scalar
multiplications. In particular, we assume that
$G\subset SL(4,\CC)$. Notation are taken from \cite[\S 8.5]{Feit-1971}
with small modifications. Here $o$ is the order of the group and
$z$ is the order of its center, $\widetilde G$ denotes some
central extension of $G$.

\begin{enumerate}
\renewcommand\labelenumi{(\Roman{enumi})}
\item
$A\times B/Z$, where $A$, $B$ are two-dimensional primitive
subgroups in $SL(2,\CC)$ and $Z$ is the central subgroup of order
$2$ which is contained in neither $A$ nor $B$,

\setcounter{enumi}{4}
\item
$\AAA_5$, $\SSS_5$,\quad $o=60,\, 120$,\ $z=1$,
\item
$\widetilde \AAA_6$, $\widetilde \SSS_6$,\quad $o=360z,\, 720z$,\ $z=2$,
\item
$\widetilde \AAA_7$,\quad $o=(\frac 12 7!) z$,\ $z=2$,
\item
$SL(2,\FF_5)$, $\widetilde \SSS_5$\quad $o=60z,\, 120z$,\ $z=2$,
\item
$SL(2,\FF_7)$,\quad $o=168z$,\ $z=2$,
\item
$G_{25920}$, \quad $o=25920z$, $z=2$,
\item
$G$ is such that $N\subset G\subset M\subset SL(4,\CC)$, where $N$
a special imprimitive subgroup of order $32$ and $M$ is its
extension by the automorphism group, $|M|=32\cdot 6!$.
\end{enumerate}
\end{emptytheorem}
Note that (XI) is the biggest class. It includes also three groups
from (I).
We discuss the rationality problem of $\PP^3/G$ for 
primitive and imprimitive subgroups in $GL(4,\CC)$
in case by case manner.

\subsection*{Notation}
Throughout this paper we will use the following notation unless
otherwise specified.
\begin{itemize}
\item
$X \bir Y$ denotes the birational equivalence of algebraic
varieties $X$ and $Y$.
\item
$\PP(a_1,\dots,a_n)$ denotes the weighted projective space.
\item
$Z(G)$ usually denotes the center of a group $G$.
\item
$\SSS_n$ and $\AAA_n$ denote the symmetric 
alternating groups on $n$ letters, respectively. Let $\SSS_n$ acts on $\CC^n$ by
permuting the coordinates $x_i$. The restriction of this
representation to the invariant hyperplane $\sum x_i=0$ we call
the \textit{standard} representation.
\item
$\TT$, $\OO$, $\II\subset SL(2,\CC)$ are binary tetrahedral,
octahedral, and icosahedral groups, respectively. There are
well-known isomorphisms 
$\TT/Z(\TT) \simeq \AAA_4$, $\OO/Z(\OO)\simeq \SSS_4$,
$\II/Z(\II) \simeq \AAA_5$,
$\TT\simeq SL(2,\FF_3)$, and $\II\simeq SL(2,\FF_5)$.
\item
If $V$ is a vector space, then $\PP(V)$ denotes its
projectivization and for an element $x\in V$, $[x]$ denotes the
corresponding point on $\PP(V)$.
\end{itemize}

For the sake of simplicity, we work over the complex number field
$\CC$. However some results can be extended for an arbitrary
algebraically closed field (at least if the characteristic is
sufficiently large).

\subsection*{Acknowledgements}
I express my thanks to V.A. Iskovskikh for many stimulating
discussions on the rationality problems.
I am grateful to F.A. Bogomolov for helpful suggestions
and  Ming-chang Kang for 
very useful corrections and comments.

\section{Groups of type (I)}
In this section we follow \cite{Kolpakov-Prokhorov-1992}. Regard
$V=\CC^4$ as the space of $2\times2$-matrices. The group
$SL(2,\CC)\times SL(2,\CC)$ acts on $V$ from the left and the right.
Therefore, for any two subgroups $G_1,\, G_2\subset SL(2,\CC)$ there is
a natural representation $\Psi\colon G_1\times G_2\to SL(V)$. Denote
$G:= \Psi(G_1, G_2)\subset SL(V)$. For $\Psi(\TT,\TT)$, $\Psi(\TT,\OO)$,
$\Psi(\TT,\II)$, $\Psi(\OO,\OO)$, $\Psi(\OO,\II)$, and
$\Psi(\II,\II)$ we get primitive groups of type (I) (they are denoted by 
$1^o$, $3^o$,
$4^o$, $5^o$, $6^o$, and $7^o$ in Blichfeldt's classification
\cite{Blichfeldt}).

\begin{theorem}[{\cite{Kolpakov-Prokhorov-1992}}]
\label{th-2-1}
The variety $\PP^3/G$ is rational for $G=\Psi (\TT,\TT)$,
$\Psi(\TT,\OO)$, $\Psi(\TT,\II)$, $\Psi(\OO,\OO)$, $\Psi(\II,\II)$.
\end{theorem}

Another proof of the rationality of $\PP^3/G$ for
$G=\Psi(\TT,\TT)$, $\Psi(\TT,\OO)$, $\Psi(\OO,\OO)$ will be given
in \S \xref{sect-Segre}.

\begin{proof}
We have
\[
\PP(V)/G \bir V/(G \cdot \CC^*) \bir SL(2,\CC)/G \bir G_1 \backslash
SL(2,\CC)/G_2.
\]
The affine variety $G_1 \backslash SL(2,\CC)$ is a homogeneous space
under the action of $SL(2,\CC)$. In all cases $G_1=\TT$, $\OO$, and
$\II$, there is a natural non-singular projective compactification
$W$ of $G_1 \backslash SL(2,\CC)$ with Picard number one. We describe
this construction below.

Denote by $M_n$ the space of binary forms of degree $n$. In the
case $G_1=\TT$, the group $\TT$ has two semi-invariants
\[
x_4,\, x_4'= t_1^4\pm 2\sqrt{-3}\ t_1^2t_2^2+t_2^4 \in M_4
\]
(see \cite[vol. 2]{Weber-1896}, \cite[\S 4.5]{Springer-1977}). Let
$W_2:=\overline{ SL(2,\CC)\cdot [x_4]}\subset \PP(M_4)=\PP^4$ be
the closure of the $SL(2,\CC)$-orbit of $[x_4]\in \PP(M_4)$. This set
is $SL(2,\CC)$-invariant and contains an open orbit isomorphic to $\TT
\backslash SL(2,\CC)$. On the other hand, there is a non-degenerate
symmetric bilinear form $q(\phantom{x},\phantom{x})$ on $M_4$ (see
\cite[vol. 2]{Weber-1896}, \cite[\S 3.1]{Springer-1977}). For
an element $a\in \TT$ of order $3$ one has $a\cdot x_4=\omega
x_4$, where $\omega$ is a primitive $3$-th root of unity.
Therefore, 
\[
q(x_4,x_4)=q(a\cdot x_4,a\cdot x_4)=\omega^2
q(x_4,x_4), 
\]
so that $q(x_4,x_4)=0$. Thus the variety $W_2$ is
defined by the equation $q(x,x)=0$. This shows that $W_2\subset
\PP^4$ is a smooth quadric.

Similarly, in cases $G_1=\OO$ and $\II$ by \cite[vol.
2]{Weber-1896} or \cite[\S 4.5]{Springer-1977}, the group $\OO$
(resp. $\II$) has a semi-invariant 
\[
x_6=t_1t_2(t_1^4-t_2^4)\in M_6,
\]

\[
\text{(resp. invariant $x_{12}=t_1t_2(t_1^{10}+11t_1^5t_2^5-t_2^{10})\in M_{12})$}. 
\]
By
\cite{Mukai-Umemura-1983} the closure $W_5:=\overline{
SL(2,\CC)\cdot [x_6]}\subset \PP(M_6)=\PP^6$ (resp.
$W_{12}:=\overline{ SL(2,\CC)\cdot [x_{12}]}\subset
\PP(M_{12})=\PP^{12}$) is a smooth compactification of
$\OO\backslash SL(2,\CC)$ (resp $\II\backslash SL(2,\CC)$).

\begin{remark}
Varieties $W_2$, $W_5$, and $W_{22}$ are smooth Fano threefolds
with $\rho=1$. One has $\Pic W=\ZZ\cdot H$ and $-K_W\sim r H$,
where $H$ is the class of hyperplane section and $r=3,\, 2$, and
$1$, respectively. The theory of $SL(2,\CC)$-varieties with an open
orbit was developed from the minimal model theoretic viewpoint in
\cite{Mukai-Umemura-1983}, \cite{Umemura-1988},
\cite{Nakano-1989b}. The main result is as follows: 
for any smooth $SL(2,\CC)$-variety $X$ with an open
orbit, there is a sequence of $SL(2,\CC)$-equivariant birational morphisms 
$X=X_1\to\cdots\to X_n=Y$, where each $X_i$ is smooth, $X_i\to X_{i+1}$
is the blow-up of a smooth $SL(2,\CC)$-invariant curve, and 
$Y$ is a so-called minimal $SL(2,\CC)$-variety.
Minimal $SL(2,\CC)$-varieties are completely classified.
In case when the stabilizer of a general point is $\II$ (resp., 
$\OO$, $\TT$), there is only one case $Y=W_{22}$ (resp., $W_5$, $W_2$).
In cases when the stabilizer is dihedral or cyclic group,
the situation is more complicated: either $Y\simeq \PP^3$ or
$Y$ has the form $\PP_S(\EEE)$, where $S$ is a smooth surface 
admitting an $SL(2,\CC)$-action 
or $\PP^1$ and $\EEE$ is an equivariant vector bundle of rank $4-\dim S$.
A complete description of pairs $(S,\EEE)$ is given.
\end{remark}

Now we consider the following cases.

\textbf{1)} $G=\Psi(\TT,\TT)$. Then 
\[
\PP(V)/\Psi(\TT,\TT)\bir \TT
\backslash SL(2,\CC)/\TT \bir W_2/\TT.
\]
The point $[x_4]\in W_2$ is
$\TT$-invariant. Projection from this point gives a birational
isomorphism $W_2/\TT \bir \PP^3/\TT$. The last variety is
obviously rational (because $\TT$ has no irreducible
four-dimensional representations).

\textbf{2)} $G=\Psi(\II,\TT)$ or $\Psi(\II,\II)$. As above,
$\PP(V)/G\bir W_{22}/G_2$ and the point $P:=[x_{12}]\in W_{22}$ is
$G_2$-invariant. Since the point $P$ is contained in an open
orbit, it does not lie on a line \cite{Mukai-Umemura-1983}. Triple
projection from $P$ (the rational map given by the linear system
$|-K_{W}-3P|$) gives a birational $G_2$-equivariant map
$W_{22}\dashrightarrow \PP^3$ \cite{Takeuchi-1989}, see also
\cite[\S 4.5]{Iskovskikh-Prokhorov-1999}. Therefore, $\PP(V)/G\bir
W_{22}/G_2 \bir \PP^3/G_2$. If $G_2=\TT$, the last variety is
rational, as in case 1). In the case $G_2=\II$, the rationality of
$\PP^3/G_2$ will be proved below.

\textbf{3)} $G=\Psi(\OO,\TT)$ or $\Psi(\OO,\OO)$. Then
$\PP(V)/G\bir W_6/G_2$ and the point $P:=[x_6]\in W_5$ is
$G_2$-invariant. Double projection from $P$ (i.e., the rational
map defined by the linear system $|-1/2 K_{W_5}-2P|$) gives a
$G_2$-equivariant rational map $W_5\dashrightarrow \PP^2$
\cite{Kolpakov-Prokhorov-1992}, \cite{Furushima-Nakayama-1989}.
There is the following diagram of $G_2$-equivariant maps
\[
\xymatrix{
W'\ar[d]^{\sigma}\ar@{-->}[r]^{\chi}&W^+\ar[d]^{\varphi}
\\
W_5\ar@{-->}[r]^{\psi}&\PP^2
}
\]
where $\sigma$ is the blow up of $P$, $\chi$ is a flop, and
$\varphi$ is a $\PP^1$-bundle. Let $S:=\sigma^{-1}(P)$ be the
exceptional divisor and let $S^+$ be its proper transform on
$W^+$. Then $S^+$ is a $G_2$-invariant section. Therefore, the
quotient $W^+/G_2\to \PP^2/G_2$ has a birational structure of
$\PP^1$-bundle. This implies that $\PP^3/G\bir W^+/G_2$ is
rational.
\end{proof}

The rationality of $\PP^3/\Psi(\OO,\II)$ is an open question.

\begin{theorem}
\label{th-O-I-st-rat}\label{cor-Psi-I}
Notation as above. 
\begin{enumerate}
\item
$\PP^3/\Psi(\OO,\II)\bir W_{22}/\OO\bir W_5/\II$.
\item
$\PP^3/\Psi(\OO,\II)$ is stably rational.
More precisely, $\PP^3/\Psi(\OO,\II)\times \PP^2$ is rational.
\end{enumerate}
\end{theorem}
\begin{proof}
(i) follows by the above arguments. We prove (ii).
Put $V:=\CC^3$. There is a faithful three-dimensional representation
$\AAA_5\to GL(V)$ that induces an action of $\AAA_5$ on $\PP^2=\PP(V)$.
We have the following fibrations
\[
\xymatrix{
&(W_5\times \PP^2)/\AAA_5 \ar[ld]_f \ar[rd]^{g} &
\\
W_5/\AAA_5&&\PP^2/\AAA_5
}
\]
where $f$ (resp. $g$) is a generically $\PP^2$ (resp. $W_5$)-bundle 
in the \'etale topology. Since the action of $\AAA_5$ on $\PP^2$ is faithful, 
the map $f$ admits a section. 
Put $X:=(W_5\times \PP^2)/\AAA_5$ and 
$K:=\CC(\PP^2/\AAA_5)$. Then $X\bir W_5/\AAA_5\times \PP^2$. On the other hand, 
$X_K:=X\otimes K$ is a smooth Fano threefold of index $2$ and degree $5$
(see \cite{Iskovskikh-Prokhorov-1999})
defined over a non-closed field $K$. 
A general pencil of hyperplane sections defines a structure of 
del Pezzo fibration of degree $5$ on $X_K$. 
By \cite[Ch. 4]{Manin-Cubic-forms-e-II} the variety $X_K$
is $K$-rational (cf. Proof of Theorem \xref{th-imprim-1-4}, case
$\Gamma_{12}^{9}$).
\end{proof}

\section{The primitive group of order $64\cdot 6!$ and the Segre cubic}
\label{sect-Segre}
In this section we prove the rationality of quotients of $\PP^3$
by primitive groups of type (XI) with two exceptions. We also give
an alternative proof of the rationality $\PP^3/G$ for groups of
type (I) for $(A,B)=(\TT,\TT)$, $(\TT,\OO)$, $(\OO,\OO)$.
This section is a modified and simplified version of 
\cite{Kolpakov-Prokhorov-1993} but basically follows 
the same idea.

\subsection*{The primitive group of order $64\cdot 6!$ (\cite{Blichfeldt})}
Let $Q_8\subset SL(2,\CC)$ be the quaternion group of order $8$. We may
assume that $Q_8=\{ \pm E, \pm I, \pm J, \pm K\}$, where
\[
I=
\begin{pmatrix}
\sqrt {-1}&0\\0&-\sqrt {-1}
\end{pmatrix},\quad
J=\begin{pmatrix} 0&\sqrt {-1}\\\sqrt {-1}&0
\end{pmatrix},\quad
K=I\cdot J=
\begin{pmatrix}
0&-1\\1&0
\end{pmatrix}.
\]
Regard $V=\CC^4$ as the space of $2\times 2$ matrices. The group
$Q_8\times Q_8$ naturally acts on $V$ by multiplications from
the left and the right. This induces a representation 
\[
\rho \colon Q_8\times Q_8\to SL(V). 
\]
The image of $\rho$ is an
imprimitive group of order $32$ isomorphic to $Q_8\times Q_8/\{\pm
1\}$. Indeed, $Q_8\times Q_8$ interchanges the one-dimensional
subspaces $V_i\subset V$ generated by $e_1:=E$, $e_2:=I$,
$e_3:=J$, and $e_4:=K$. Denote by $N$ the subgroup of order $64$
in $SL(V)$ generated by $\rho(Q_8\times Q_8)$ and scalar
multiplications by $\sqrt{-1}$. Let $M\subset SL(V)$ be the
normalizer of $N$. We study the rationality question for different
subgroups of $M$.

The group $M$ naturally acts on $\wedge^2V\simeq \CC^6$ as an
imprimitive group (see \cite{Blichfeldt}). This easily follows
from the fact that the image of $N$ in $GL(\Lambda^2V)$ is an
Abelian group. The eigenvectors
\[
\begin{array}{lll}
w_{12}^+:=e_1\wedge e_2 + e_3\wedge e_4,\ & w_{13}^+:=e_1\wedge
e_3 + e_2\wedge e_4,\ & w_{14}^+:=e_1\wedge e_4 + e_2\wedge e_3,
\\
w_{12}^-:=e_1\wedge e_2 - e_3\wedge e_4,& w_{13}^-:=e_1\wedge e_3
- e_2\wedge e_4, & w_{14}^-:=e_1\wedge e_4 - e_2\wedge e_3
\end{array}
\]
form a basis of $\wedge^2 V$. This gives a decomposition
\begin{equation}
\label{eq-Segre-exseq-gr-0}
\bigwedge\nolimits^2V=\bigoplus_{i=1}^6 W_i, \qquad W_i=\CC \cdot w_{j,k}^{\pm} 
\end{equation}
such that $N\cdot W_i=W_i$. The group $M/N$ permutes subspaces
$W_i$ and is isomorphic to the symmetric group $\mathfrak{S}_6$.
Thus we have the exact sequence
\begin{equation}
\label{eq-Segre-exseq-gr}
1\longrightarrow N \longrightarrow M
\stackrel{\pi}{\longrightarrow} \mathfrak{S}_6\longrightarrow 1
\end{equation}
and any finite group containing $N$ as a normal subgroup is
uniquely determined by its image under $\pi\colon M\to
\mathfrak{S}_6$.

\begin{example}
\label{exam-Segre}
It is easy to see that the matrix
$S:=\frac{1+\sqrt{-1}}{\sqrt2}\Diag (\sqrt{-1},\sqrt{-1},1,1)$ is
contained in $M$. The action on $\wedge^2 V$ is as follows
\[
\begin{array}{lll}
w_{12}^+ \rightarrow -\sqrt{-1} \, w_{12}^-,\quad & w_{13}^+
\rightarrow -w_{13}^+,\quad & w_{14}^+\rightarrow -w_{14}^+,
\\
w_{12}^- \rightarrow -\sqrt{-1} \, w_{12}^+,& w_{13}^-\rightarrow
-w_{13}^-,& w_{14}^-\rightarrow -w_{14}^-.
\end{array}
\]
Therefore, $S$ is a transposition $(1,2)$ in $M/N=\SSS_6$.
Similarly, the matrix $B:=\frac{1+\sqrt{-1}}{\sqrt2}\Diag
(1,1,1,-1)$ also normalizes $N$. The corresponding substitution in
$\SSS_6$ is $(1,2)(3,4)(5,6)$ (odd element of order $2$).
\end{example}

\begin{remark}
The group $N/Z(N)=Q_8\times Q_8/\{\pm 1,\pm 1\}$ is Abelian of
order $16$ isomorphic to $(\FF_2)^4$. The exact sequence
\eqref{eq-Segre-exseq-gr} induces an embedding
$\SSS_6\hookrightarrow \Aut \FF_2^4\simeq SL(4,\FF_2)$.
\end{remark}

\subsection*{The results}
\begin{theorem}[\cite{Kolpakov-Prokhorov-1993}]
\label{cor-Igusa-quotient}
The quotient $\PP^3/N$ is isomorphic to the Igusa quartic
$\mathcal I_4\subset \PP^4$ \textup(see below\textup).
\end{theorem}

\begin{theorem}[\cite{Kolpakov-Prokhorov-1993}]
\label{th-Segre-1}
Let $G\subset SL(4,\CC)$ be a finite group having a normal
imprimitive subgroup $N$ of order $64$. Then $\PP^3/G$ is
birationally isomorphic to the quotient of the Segre cubic
$\mathcal S_3\subset \PP^4$ by $G/N$, where the action of $G/N$ on
$\PP^4$ is given by the composition of the canonical embedding
$\pi\colon G/N\hookrightarrow \SSS_6$, an outer automorphism
$\lambda\colon \SSS_6\to \SSS_6$ \textup(see, e.g.,
\cite{Miller-1958}\textup), and the standard action of $\SSS_6$ on
$\mathcal S_3$.
\end{theorem}

\begin{corollary}[\cite{Kolpakov-Prokhorov-1993}]
\label{th-Segre-2}
In the above notation, the variety $\PP^3/G$ is rational except
possibly for the following two groups \textup(we use Blichfeldt's
notation \cite{Blichfeldt}\textup):
\begin{enumerate}
\item[$20^o$]
$|G|=64\cdot 360$ and $G/N\simeq \AAA_6$,
\item[$17^o$]
$|G|=64\cdot 60$, $G/N\simeq \AAA_5$, and $\pi(G/N) \subset\SSS_6$
as a transitive subgroup.
\end{enumerate}
\end{corollary}

\begin{corollary}
\label{th-Segre-3}
In cases $20^o$ and $17^o$ above, the variety $\PP^3/G$ is
birationally isomorphic to $\mathcal S_3/\lambda \comp \pi(G/N)$.
\end{corollary}
Note that in case $17^o$ one of the $x_i$ is $\lambda \comp \pi(G/N)$-invariant.

\subsection*{Segre cubic}
Regard the variety $\mathcal S_3$ given by the following equations
\begin{equation}
\label{eq-Segre-cubic}
\begin{array}{lll}
x_1+x_2+x_3+x_4+x_5+x_6&=0,
\\[5pt]
x_1^3+x_2^3+x_3^3+x_4^3+x_5^3+x_6^3&=0.
\end{array}
\end{equation}
as a cubic hypersurface in $\PP^4$. This cubic satisfies many
remarkable properties (see \cite{Segre-1963}, \cite[Ch.
8]{Semple-Roth-1949}, \cite{Dolgachev-Ortland-1989}) and is called
the \emph{Segre cubic}. For example, any cubic hypersurface in
$\PP^4$ has at most ten isolated singular points, this bound is
sharp and achieved exactly for the Segre cubic (up to projective
isomorphism). The symmetric group $\SSS_6$ acts in $\mathcal S_3$
in the standard way. Moreover, it is easy to show that $\Aut(\mathcal
S_3)=\SSS_6$ (see e.g. \cite{Finkelnberg-1987}). We refer to
\cite{Dolgachev-Ortland-1989} and \cite{Koike} for further
interesting properties of $\mathcal S_3$.

The singular locus of $\mathcal S_3$ consists of ten nodes given
by
\begin{equation}
\label{eq-Segre-cubic-sing}
x_{i_1}=x_{i_2}=x_{i_3}=-x_{i_4}=-x_{i_5}=-x_{i_6},
\end{equation}
where $\{i_1,i_2,i_3,i_4,i_5,i_6\}=\{1,2,3,4,5,6\}$. We denote
such a point by $\point{i_1,i_2,i_3}{i_4,i_5,i_6}$. For example,
$\point{123}{345}=(1,1,1,-1,-1,-1)$. It is easy to see that
$\point{i_1,i_2,i_3}{i_4,i_5,i_6}=\point{j_1,j_2,j_3}{j_4,j_5,j_6}$
if and only if corresponding matrices are obtained from each other
by permutations of rows and elements in each row. Hence there is an
$\SSS_6$-equivariant $1$-$1$-correspondence
\[
\Sing(\mathcal S_3) \longleftrightarrow \{\text{Sylow $3$-groups
in $\SSS_6$}\}.
\]
Since the action of $\SSS_6$ on $\Sing(\SSS_3)$ is transitive, the
stabilizer $\St(\point{i_1,i_2,i_3}{i_4,i_5,i_6})$ of
$\point{i_1,i_2,i_3}{i_4,i_5,i_6}$ is a group of order $72$. For
example, $\St(\point{123}{345})$ is generated by $\SSS_3\times
\SSS_3$ and $(1,3)(2,4)(3,5)$.

Further, there are $15$ planes on $\mathcal S_3$ (see e.g.
\cite{Finkelnberg-1987}). Each of them is given by equations
\begin{equation}
\label{eq-Segre-cubic-planes}
x_{i_1}+x_{i_4}=x_{i_2}+x_{i_5}=x_{i_3}+x_{i_6}=0,
\end{equation}
where $\{i_1,i_2,i_3,i_4,i_5,i_6\}=\{1,2,3,4,5,6\}$. Denote such a
plane by $\plane{i_1,i_2,i_3}{i_4,i_5,i_6}$ and the set of all the
planes on $\mathcal S_3$ by $\Omega$. It is easy to see that
$\plane{i_1,i_2,i_3}{i_4,i_5,i_6}=\plane{j_1,j_2,j_3}{j_4,j_5,j_6}$
if and only if corresponding matrices are obtained from each other
by permutations of columns and elements in each column. Thus there
is an $\SSS_6$-equivariant $1$-$1$-correspondence
\[
\Omega \longleftrightarrow \{(i_1,i_4)(i_2,i_5)(i_3,i_6)\}\subset
\SSS_6.
\]
For $\sigma=(i_1,i_4)(i_2,i_5)(i_3,i_6)$, we often will write
$\Pi(\sigma)$ instead of $\plane{i_1,i_2,i_3}{i_4,i_5,i_6}$.

\begin{corollary}
$\St(\Pi(\sigma))=Z(\sigma)$ and $|Z(\sigma)|=48$.
\end{corollary}

Comparing \eqref{eq-Segre-cubic-planes} and
\eqref{eq-Segre-cubic-sing} one can see the following facts.
\begin{enumerate}
\item[a)]
Every plane $\Pi(\sigma)$ contains exactly four singular points
and there are six planes passing trough every singular point.
\item[b)]
$\Pi(\sigma_1)\cap\Pi(\sigma_2) =\left\{
\point{i_1,i_2,i_3}{i_4,i_5,i_6} \right\}$ $\Longleftrightarrow$
$\sigma_1\comp \sigma_2=(i_1,i_2,i_3)(i_4,i_5,i_6)$.
\item[c)]
$\Pi(\sigma_1)\cap\Pi(\sigma_2)$ is a line $\Longleftrightarrow$
$\sigma_1$ and $\sigma_2$ contain a common transposition.
\end{enumerate}
Recall that there exists an outer automorphism $\lambda\colon
\SSS_6\to \SSS_6$ (see, e.g., \cite{Miller-1958}). For any
standard embedding $\SSS_5 \subset \SSS_6$, the subgroup
$\lambda(\SSS_5)$ is a ``non-standard'' transitive subgroup
isomorphic to $\SSS_5$. We need the following simple lemmas.

\begin{lemma}[cf. {\cite[pp. 169--170]{Semple-Roth-1949}}]
\label{lemma-tran--S5}
Let $G\subset \SSS_6$ be a subgroup of order $120$
\textup(isomorphic to $\SSS_5$\textup). One of the following
holds:
\begin{enumerate}
\item
$G$ is not transitive, then $G\cap \St(\Pi(\sigma))$ is of order
$8$ and $G$ acts on $\Omega$ transitively;
\item
$G$ is transitive, then the order of $G\cap \St(\Pi(\sigma))$ is
either $12$ or $24$ and $\Omega$ splits into two $G$-orbits
$\Omega'$ and $\Omega''$ consisting of $10$ and $5$ planes,
respectively. A plane $\Pi(\sigma)$ is contained in $\Omega'$ if
and only if $\sigma\in G$. Moreover, every two planes from
$\Omega''$ intersect each other at a \textup(singular\textup)
point.
\end{enumerate}
\end{lemma}
\begin{proof}
(ii) Assume that there are two planes $\Pi(\sigma_1),
\Pi(\sigma_2)\in \Omega''$ such that
$\Pi(\sigma_1)\cap\Pi(\sigma_2)$ is a line. Then $\sigma_1$ and
$\sigma_2$ contain a common transposition and $\sigma_1,
\sigma_2\notin G$. This is equivalent to that
$\sigma_1\comp\sigma_2$ is an element of order $2$. But then
$\lambda(\sigma_1)$, $\lambda(\sigma_2)$ are two transpositions
such that $\lambda(\sigma_1)\comp\lambda(\sigma_2)$ is an element
of order $2$ and both $\lambda(\sigma_1),\, \lambda(\sigma_2)$ are
not contained in a standard subgroup $\SSS_5\subset \SSS_6$.
Clearly, this is impossible.
\end{proof}

\begin{lemma}
\label{lemma-Pi-pencil-Q}
Let $\Pi=\Pi(\sigma)$ be a plane on $\mathcal S_3$ and let
$\mathcal H$ be the pencil of hyperplane sections through $\Pi$.
Let $\mathcal Q$ be the pencil of residue quadrics to $\Pi$
\textup(i.e., $\mathcal H=\Pi+\mathcal Q$\textup). Then the base
locus of $\mathcal Q$ consists of four singular points of $\mathcal
S_3$ contained in $\Pi$ and a general member of $\mathcal Q$ is
smooth.
\end{lemma}

\begin{proof}
It is sufficient to check statement only for one plane. In this
case it can be done explicitly.
\end{proof}

\subsection*{Igusa quartic}
Consider the dual map $\Psi\colon \mathcal S_3\dashrightarrow
\PP^{4*}$ sending a smooth point $P\in \mathcal S_3$ to the
tangent space $T_{P,\mathcal S_3}$. Let $\mathcal
S_3^*=\Psi(\mathcal S_3)$ be the \emph{dual variety}.

\begin{lemma}[\cite{Koike}, \cite{Dolgachev-Ortland-1989}]
\label{lemma-Igusa-bir}
The map $\Psi\colon \mathcal S_3 \dashrightarrow \mathcal S_3^*$
is birational and $\deg \mathcal S_3^*=4$.
\end{lemma}
\begin{proof}
Assume that $\Psi$ is not birational. Then for a general $P\in
\mathcal S_3$ there is at least one point $P'\neq P$ such that
$\Psi(P)=\Psi(P')$. This means that tangent spaces $T_{P,\mathcal
S_3}$ and $T_{P',\mathcal S_3}$ coincide. We may assume that $P,\,
P'$ are not contained in any plane $\Pi(\sigma)$. The line $L$
passing through $P$ and $P'$ is contained in $\mathcal S_3$. This
line meets some plane $\Pi=\Pi(\sigma)\subset\mathcal S_3$. Thus
the linear span $\langle L,\, \Pi\rangle$ is a hyperplane in $\PP^4$. It is easy
to see that $\langle L,\, \Pi\rangle\cap \mathcal S_3=\Pi\cup Q$,
where $Q$ is a two-dimensional quadric. Clearly,
$T_{P,Q}=T_{P',Q}$. This is possible only if $Q$ is a quadratic
cone. By Bertini's theorem its vertex is contained in $\Pi$. This
contradicts Lemma \xref{lemma-Pi-pencil-Q}.

Therefore, $\Psi$ is birational and a general tangent hyperplane
section $T_{P,\mathcal S_3}\cap \mathcal S_3$ has exactly one
node. This implies that a general pencil
$\mathcal H$ of hyperplane sections of $\mathcal S_3$ is a
Lefschetz pencil. Let $H_1,\dots,H_{r}$ be
singular members of $\mathcal H$. Ten of them pass through singular
points of $\mathcal S_3$. Hence the degree of $\mathcal
S_3^*\subset \PP^{4*}$ is equal to $r-10$. Finally, it is easy to
compute that the topological Euler number of a cubic hypersurface
in $\PP^4$ having only nodes as singularities is equal to
$m-6$, where $m$ is the number of singular points. Using this fact
one can see that $\chit(\mathcal S_3)=4= 9(2-r)+8r$, so $r=14$ and
$\deg \mathcal S_3^*=4$.
\end{proof}
Therefore, $\mathcal S_3^*$ is a hypersurface of degree $4$
in $\PP^4$.
This famous quartic is called the \emph{Igusa quartic} (cf.
\cite{Dolgachev-Ortland-1989}, \cite{Koike}) and denoted by
$\mathcal I_4$. An interesting fact is that $\mathcal I_4$ is the
compact moduli space of Abelian surfaces with the level two
structure. By the above, there is an $\SSS_6$-equivariant
birational map $\Psi\colon \mathcal S_3 \dashrightarrow \mathcal
I_4$. We need the following characterization of $\mathcal I_4$ in
terms of the action of $\SSS_6$.

\begin{lemma}
\label{lemma-igusa-planes-lines}
The image of every plane $\Pi=\Pi(\sigma)$ is a line on $\mathcal
S_3^*$. Therefore, $\mathcal S_3^*$ contains an $\SSS_6$-invariant
configuration of $15$ lines such that the action on these lines is
transitive. Moreover, the above $15$ lines form the singular locus
of $\mathcal S_3^*$ and the stabilizer of such a line is
$Z(\sigma)$, where $\sigma=(i_1,i_2)(i_3,i_4)(i_5,i_6)$.
\end{lemma}
\begin{proof}
The tangent space $T_{P,\mathcal S_3}$ at $P\in \Pi$ contains
$\Pi$. Hence $T_{P,\mathcal S_3}$ is contained in the pencil of
hyperplane sections passing through $\Pi$. The fact that
$\Psi(\Pi)\subset \Sing(\mathcal S_3^*)$ can be checked by direct
computations, see \cite{Koike}.
\end{proof}

\begin{lemma}
\label{lemma-igusa-planes-lines-1}
Let $X\subset \PP^4$ be an $\SSS_6$-invariant quartic under the
standard action of $\SSS_6$ on $\PP^4$. Assume that there is an
odd element $\sigma\in \SSS_6$ of order $2$ such that
\begin{enumerate}
\item
the fixed point locus of $\sigma$ on $\PP^4$ is a disjointed union
of a line $L$ and a plane $\Pi$,
\item
$L$ is contained in $X$.
\end{enumerate}
Then $X$ is the Igusa quartic \textup(up to projective
equivalence\textup) and $L\subset \Sing(X)$.

\end{lemma}
\begin{proof}
It is clear that $X$ must be defined by an invariant $\psi$ of
degree $4$. On the other hand, modulo $\sum x_i$ there are exactly
two linearly independent invariants of degree $4$: $s_4$ and
$s_2^2$, where $s_d:=\sum x_i^d$. Hence, $\psi=\alpha s_4+\beta
s_2^2$, $\alpha,\beta\in \CC$.

Applying a suitable automorphism of $\SSS_6$ we may assume that
$\sigma=(12)(34)(56)$. Therefore, $L$ is given by the following
equations:
\[
x_1=x_2,\quad x_3=x_4,\quad x_5=x_6,\quad x_1+x_3+x_5=0.
\]
On the other hand, in the pencil $\psi=\alpha s_4+\beta s_2^2$
there is exactly one quartic containing such an $L$.
\end{proof}

\subsection*{Invariants of $N$}

\begin{lemma}
The ring of invariants $\CC[x,y,z,u]^{N}$ is generated by the
following elements of degree $4$:
\[
\begin{array}{ll}
f_1=x^2u^2+y^2z^2,\quad &f_2=x^2z^2+y^2u^2, \quad
f_3=x^2y^2+z^2u^2,
\\[8pt]
f_4=x^4+y^4+z^4+u^4,\quad &f_5=xyzu.
\end{array}
\]
The only relation is
\begin{multline}
\label{eq-Igusa-quartic-1}
\left(f_4+2f_1+2f_2+2f_3\right)\Bigl(f_1f_2f_3-f_4f^2_5+2f^2_5(f_1+f_2+f_3)\Bigr)-
\\
-\left(f_1f_2+f_2f_3+f_3f_1+4f^2_5\right)^2=0.
\end{multline}
\end{lemma}
\begin{proof}
Indeed, it is easy to see that all the $f_i$ are $N$-invariants.
Hence, $\CC[f_1,\dots,f_5]\subset \CC[x,y,z,t]^{N}$. Now it is
sufficient to show that
\[
\left[\CC(x,\, y,\, z,\, u):\CC(f_1,\, f_2,\, f_3,\, f_4,\,
f_5)\right]=|N|=64.
\]
Using the following relations
\[
\left[\CC(f_1,\, f_2,\, f_3,\, f_4,\, f_5):\CC(f_1,\, f_2,\,
f_3,\, f_4,\, f_5^2)\right]=2,
\]
\[
\left[\CC(f_1,f_2,f_3,f_4,f_5^2) :
\CC(xu+yz,xz+yu,xy+zu,x^2+y^2+z^2+u^2,f_5)\right] =16,
\]
\begin{multline*}
[\CC(xu+yz,\, xz+yu,\, xy+zu,\, x^2+y^2+z^2+u^2,\, f_5) :
\\
: \CC(xu+yz,\, xz+yu,\, xy+zu,\, x+y+z+u,\, f_5)]=2
\end{multline*}
we get the desired equality. Relation \eqref{eq-Igusa-quartic-1}
can be checked directly. It is easy to see that the quartic given
by \eqref{eq-Igusa-quartic-1} is smooth in codimension one. Hence
it is irreducible and \eqref{eq-Igusa-quartic-1} is the only
relation between the $f_i$.
\end{proof}

\begin{lemma}
The quartic $X:=\PP^3/N$ given by \eqref{eq-Igusa-quartic-1} is
singular along $15$ lines $L_1,\dots,L_{15}$. Moreover, $\SSS_6$
acts on the $L_i$ transitively.
\end{lemma}

\begin{proof}
This can be obtained immediately from equation
\eqref{eq-Igusa-quartic-1} but we prefer to use the quotient structure 
on $X$. The group $N/Z(N)$ is an Abelian group of
order $16$ isomorphic to $(\muu_2)^4$. Every non-trivial element
$a\in N/Z(N)$ fixes points on two lines $\bar L_a^{i}\subset \PP^3$, $i=1,\,
2$. Thus in $\PP^3$ there are $30$ lines whose general points
have stabilizer of order $2$. The images of the $\bar L_a^i$ are
contained in the singular locus of $X$.
\end{proof}

\begin{proof}[Proof of Theorem \xref{cor-Igusa-quotient}]
Put $X:=\PP^3/N$ and fix the embedding $X\hookrightarrow \PP^4$ by
the $f_i$'s. Clearly, the group $\SSS_6\simeq M/N$ naturally acts on
$X$. Since $\Pic X\simeq \ZZ\cdot (-K_X)$, the action is induced
by a linear action on $\PP^4$. The transposition $S$ from Example
\xref{exam-Segre} acts on the $f_i$ by the diagonal matrix $\Diag
(1,1,-1,-1,1)$. In particular, it is not a reflection. This shows
that up to scalar multiplication the action of $\SSS_6$ on
$X\subset \PP^4$ is given by the composition of the standard
action of $\SSS_6$ on $\PP^4$ and the outer automorphism
$\lambda$. Now it is sufficient to show that $\Sing(X)$ contains a
line $L$ such as in Lemma \xref{lemma-igusa-planes-lines-1}.
Indeed, the fixed point locus of $S$ on $\PP^4$ consists of
disjointed union of the plane $f_3=f_4=0$ and the line
$f_1=f_2=f_5=0$ contained in $X$. Moreover, $X$ is singular along
$f_3=f_4=0$. Thus Lemma \xref{lemma-igusa-planes-lines-1} can be
applied and $X\simeq \mathcal S_3$.
\end{proof}

\begin{proof}[Proof of Theorem \xref{th-Segre-1} and Corollary \xref{th-Segre-2}]
Consider any subgroup $G$ such that $N\subset G\subset M$. By Theorem
\xref{cor-Igusa-quotient} the quotient $\PP^3/N$ is the Igusa
quartic $\mathcal I_4\subset \PP^4$, where the action of
$G/N=\SSS_6$ on $\PP^4$ is given by the composition of the
canonical embedding $\pi\colon G/N\hookrightarrow M/N=\SSS_6$, an
outer automorphism $\lambda\colon \SSS_6\to \SSS_6$, and the
standard action of $\SSS_6$ on $\mathcal \PP^4$. Thus $\PP^3/G
\bir \mathcal I_4/ (G/N)$. By Lemma \xref{lemma-Igusa-bir} there
is an $\SSS_6$-equivariant birational map $\mathcal I_4
\dashrightarrow \mathcal S_3$. Hence $\PP^3/G \bir \mathcal S_3/
(G/N)$. This proves Theorem \xref{th-Segre-1}.

To prove Corollary \xref{th-Segre-2} we consider two cases.

\begin{emptytheorem}
\label{a)} 
\textbf{$G/N$ has a fixed point $P\in \Sing (\mathcal S_3)$.} Projection
from this point is $G/N$-equivariant, and therefore $\mathcal S_3/
(G/N)$ is birationally equivalent to $\PP^3/(G/N)$, where $G/N$ is
a group of order $\le 72$. Corollary
\xref{th-Segre-2} in this case follows from the rationality
of $\PP^3/(G/N)$.
\end{emptytheorem}

\begin{emptytheorem}
\label{b)} 
\textbf{$\pi(G)$ is a subgroup of some $\SSS_5\subset \SSS_6$, a
standard (non-transitive) permutation group.} This exactly means
that one of the $W_i$ (see \eqref{eq-Segre-exseq-gr-0}) is an
eigenspace for $G$. Then $\PP^3/G \bir \mathcal S_3/ (G/N)$, where
$G/N$ is embedded into a transitive $\SSS_5$. By Lemma
\xref{lemma-tran--S5} there is a $G$-invariant set $\Omega''$ of
$5$ planes $\Pi_i=\Pi(\sigma_i)$ and $\Pi_i\cap \Pi_j\subset
\Sing(\mathcal S_3)$ for $i\neq j$. We claim that there are
exactly two planes $\Pi_i\in \Omega''$ passing through every
singular point. Indeed, this follows by the fact that $\SSS_5$
transitively acts on the set of Sylow $3$-subgroups in $\SSS_6$.
Further, for every singular point $P\in \mathcal S_3$, there is
exactly two small (possibly non-projective) resolutions
$\mu_P\colon \bar {\mathcal S}_3\to \mathcal S_3$ and
$\mu_P'\colon \bar {\mathcal S}_3'\to \mathcal S_3$. Recall that a
birational contraction is said to be \textit{small} if it does not
contract any divisors. For one of them, say for $\mu_P$, the
proper transforms $\bar\Pi_i$ of planes $\Pi_i\in \Omega''$ does
not meet each other over $\mu_P^{-1}(P)$. Thus there is exactly
one small resolution $\mu\colon \bar {\mathcal S}_3\to \mathcal
S_3$ of all singular points such that the proper transforms $\bar
\Pi_i$ of the planes $\Pi_i\in \Omega''$ are disjointed. This
resolution must be $\SSS_5$-equivariant (cf. \cite[\S
5]{Finkelnberg-1987}). By our construction, $\bar \Pi_i\simeq
\PP^2$. It is easy to check that the $\bar \Pi_i$ satisfy the
contractibility criterion, i.e., $\bar \Pi_i|_{\bar \Pi_i}\simeq
\OOO_{\PP^2}(-1)$. So there is an $\SSS_5$-equivariant contraction
$\varphi\colon \bar {\mathcal S}_3\to W$ of Moishezon complex
manifolds. Since $\Pic \bar {\mathcal S}_3\simeq \ZZ^{\oplus 6}$ (see
\cite{Finkelnberg-1987}), we have $\Pic W\simeq \ZZ$. Obviously,
the anticanonical divisor $-K_W$ is
effective and divisible by four in $\Pic W$. By
\cite{Nakamura-1987} the variety $W$ is projective and $W\simeq
\PP^3$. Therefore, we have $\PP^3/G \bir \mathcal S_3/(G/N)\bir
\PP^3/(G/N)$, where $(G/N)$ is a subgroup of $\SSS_5$. The latter
reduces the question of rationality of $\PP^3/G$ to a smaller
group $G$ of order $\le 120$. This will be discussed below.

In fact, it will be shown in \S \xref{sect-simple-and-other} that
there is another (different from $\mu\colon \bar {\mathcal S}_3\to
\mathcal S_3$) small projective $G/N$-equivariant resolution
$\mu^+\colon \bar {\mathcal S}_3^+\to \mathcal S_3$ and there is a
$G/N$-equivariant $\PP^1$-bundle structure on $\bar {\mathcal
S}_3^+$ (see Proof of Proposition \xref{th-A5}).
\end{emptytheorem}

To finish the proof of Corollary \xref{th-Segre-2} we note that, except for
$20^o$ and $17^o$, there are only two groups which do not satisfy
conditions \xref{a)} or \xref{b)} above: $\pi (G)=\SSS_6$ and $\pi(G)\subset
\SSS_6$ is a transitive $\SSS_5$. The ring of invariants
$\CC[\mathcal S_3]^{\SSS_6}$ is generated by symmetric functions
$s_2$, $s_4$, $s_5$, $s_6$, where $s_k=\sum_{i=1}^6 x_i^k$.
Therefore, $\mathcal S_3/\SSS_6\simeq \PP(2,4,5,6)$ is rational.
Consider the case when $G/N\simeq \SSS_5$ and $\pi(G)\subset
\SSS_6$ is a transitive $\SSS_5$. Then $G/N=\SSS_5$ fixes some of
$x_1,\dots,x_6$ on $\mathcal S_3$. Assume that $\SSS_5\cdot
x_6=x_6$. The ring of invariants $\CC[\mathcal S_3]^{\SSS_6}$ is
generated by $x_6$ and $s_1'$, \dots, $s_5'$, where
$s_k'=\sum_{i=1}^5 x_i^k$. Thus, $\CC[\mathcal S_3]^{\SSS_6}\simeq
\CC[s_2',s_4',s_5',x_6]$ and $\mathcal S_3/\SSS_5\simeq
\PP(2,4,5,1)$ is rational. This proves Corollary
\xref{th-Segre-2}.
\end{proof}

\begin{remark}
Another approach to the treatment of case \xref{b)} is to 
note that there is a $G/N$-equivariant $\PP^1$-bundle structure 
on $\mathcal S_3$.
Indeed, consider the family of lines 
$\mathcal L=\mathcal L(\mathcal S_3)$ on $\mathcal S_3$.
Let $\mathcal L_1,\dots,\mathcal L_r$ be all covering
irreducible components (i.e., components $\mathcal L_i$ such that 
there is a line from $\mathcal L_i$ passing through a general 
point of $\mathcal S_3$). Since there is at most $6$ 
lines passing through a general point, $r\le 6$. On the other hand, 
there is a hyperplane section of the form $\Pi_1+\Pi_2+\Pi_3$,
where $\Pi_i$ are planes. A general line from a covering family 
meets exactly one of the planes $\Pi_1$, $\Pi_2$, $\Pi_3$. This shows 
that $r\ge 3$. Using the action of $\SSS_6$ of $\mathcal S_3$ one can 
see that $r=6$. Therefore, each family $\mathcal L_i$ gives us 
(birationally) a $\PP^1$-bundle structure on $\mathcal S_3$.
Now it remains to note that 
if $G$ is such as in \xref{b)}, then one of families $\mathcal L_i$
is $G/N$-invariant. 
\end{remark}

\section{Groups of types (V)-(X)}
\label{sect-simple-and-other}
\subsection*{Case (X)}
Recall that an element $g\in GL(n,\CC)$ of finite order is said to
be (complex) \textit{reflection} if exactly $n-1$ eigenvalues are
equal to $1$. For convenience of the reader we 
recall the following well-known theorem of
Chevalley and Shephard-Todd:

\begin{theorem}[see {\cite{Shephard-Todd}}, 
{\cite{Chevalley-1955}}, see also \cite{Springer-1977}]
\label{th-Chevalley}
Let $V=\CC^n$ and let $G\subset GL(V)$ be a finite subgroup. The following are
equivalent:
\begin{enumerate}
\item
$G$ is generated by reflections,
\item
$V/G\simeq \CC^n$.
\end{enumerate}
\end{theorem}

\begin{proof}[Outline of Proof]
Let $W:=V/G$ and $f\colon V\to W$ be the quotient morphism.

(ii) $\Longrightarrow$ (i). Assume the converse. Let $G_0\subset G$ be the maximal 
subgroup generated by reflections. Then the morphism $V/G_0\to W$ 
is \'etale over $W\setminus Z$, where $Z$ is of codimension of least two.
On the other hand, $\pi_1(V\setminus Z)=\{1\}$, a contradiction.

(i) $\Longrightarrow$ (ii). 
Put $R:=\CC[x_1,\dots , x_n]$. First we claim that 
$R$ is a free $R^G$-module. Let $I\subset R$ be 
the ideal generated by homogeneous invariants of positive degree.

\begin{lemma}
Assume that for some homogeneous elements
$y_i\in R$ and $z_i\in R^G$ the following relation holds
\begin{equation}
\label{eq-lemma-Che}
z_1y_1+\cdots+z_my_m=0.
\end{equation}
If $z_1\notin R^Gz_2+\cdots+R^Gz_m$, then $y_1\in I$.
\end{lemma}

\begin{proof}
We can the take $y_i$ so that $\deg y_1$ 
is minimal. If $\deg y_1=0$, then applying the map 
$\frac1{|G|}\sum_{g\in G} g$ to \eqref{eq-lemma-Che}
we obtain $z_1\in R^Gz_2+\cdots+R^Gz_m$. 
Assume that $\deg y_1>0$.
Let $s\in G$ be a reflection and let $\{l_s=0\}$ be its 
fixed hyperplane. 
For any homogeneous $f\in R$, the polynomial $s\cdot f-f$
is divisible by $l_s$. 
Define a map
\[
\Delta_s\colon R_d\to R_{d-1}, \qquad \Delta_s(f)=(s\cdot f-f)/l_s.
\]
It is easy to check that 
\[
\Delta_s(f_1f_2)=(\Delta_sf_1)f_2,\quad \forall f_1\in R,\quad 
\forall f_2\in R^G.
\]
This gives us 
\[
z_1\Delta_s(y_1)+\cdots+z_m\Delta_s(y_m)=0.
\]
Since $\deg \Delta_s(y_i)<\deg y_i$, we may assume that
$\Delta_s(y_1)\in I$. Therefore, $s\cdot y_1-y_1\in I$
for all reflections $s\in G$. This implies that 
$g\cdot y_1-y_1\in I$
for all $g\in G$, so $y_1\in I$. 
\end{proof}

Take homogeneous elements $y_i\in R$ so that the images $\bar y_i$ in $R/I$ form 
a basis over $\CC$. It is clear that the $y_i$ generate 
$R$ as an $R^G$-module. By the above lemma these $y_i$
are linearly independent over $R^G$. Indeed, 
if $z_1y_1+\cdots+z_my_m=0$ for some $z_i\in R^G$, then
$z_1=z_2u_2+\cdots + z_mu_m$, $u_i\in R^G$. So, 
\[
z_2(y_2+y_1u_2)+\cdots+z_m(y_m+y_1u_m)=0
\] 
and we can apply the induction by $m$.
Since $R$ is integral over $R^G$, the basis $y_i$ is finite.
This proves our claim.

Further, let $J$ be the 
the ideal of $R^G$ generated by homogeneous elements of positive degree.
Since $R^G$ is a Noetherian algebra, $J$ is finitely generated.
Take a minimal system of generators $f_1,\dots, f_r$.
It is clear that $f_1,\dots, f_r$ generate $R^G$ as $\CC$-algebra.
On the other hand, one can check that they are algebraically 
independent.
This proves theorem.
\end{proof}
Note that the statement of Theorem \xref{th-Chevalley} fails
if the characteristic of the base field divide $|G|$. 
However the field of invariants is rational in this case 
under the additional assumption that $G$ is irreducible \cite{Kemper-1999}.

\begin{corollary}
Let $G\subset GL(n,\CC)$ be a finite group generated by reflections. Then 
$\CC[x_1,\dots,x_n]^G\simeq \CC[f_1,\dots,f_n]$, where $f_1,\dots,
f_n$ are homogeneous polynomials.
\end{corollary}
Numbers $d_i:=\deg f_i$ are called the \textit{degrees} of $G$.
They are uniquely determined by $G$. Thus for any group generated
by reflections we have $\CC^n/G\simeq \CC^n$ and $\PP^{n-1}/G$ is
a weighted projective space $\PP(d_1,\dots,d_n)$.

The list of all complex reflection groups can be found in 
\cite{Shephard-Todd}, \cite{Cohen-1976}. According to this list
there is a finite group No. 32 of order $25920\cdot 6$ generated by complex
reflections of order $3$. The degrees of this group are $12$,
$18$, $24$, $30$ and the intersection with
$SL(4,\CC)$ is exactly our first group.
Therefore, the quotient $\PP^3/G_{25920}\simeq \PP(12, 18,24, 30)
\simeq \PP(2, 3,4, 5)$ is rational.

\subsection*{The group $SL(2,\FF_5)$}
Let $\delta\colon SL(2,\FF_5) \hookrightarrow SL(2,\CC)$ be a
faithful representation whose image is the icosahedron group
$\II$. For short, we identify $SL(2,\FF_5)$ with $\II$. Then
$S^3 \delta \colon \II \to SL(4,\CC)$ is an irreducible faithful
representation as in (VIII). This gives as the action of $\AAA_5=\II /\{\pm
E\}$ on $\PP^3$ which leaves a rational cubic curve $C$ invariant.
Since $\II$ has a faithful two-dimensional representation,
$\PP^3/\II$ is stably rational. We prove more:

\begin{theorem}[\cite{Kolpakov-Prokhorov-1987}, \cite{Kolpakov-Prokhorov-1992}]
\label{th-SL25}
$\PP^3/\II$ is rational.
\end{theorem}
\begin{proof}
Let $\sigma \colon X\to \PP^3$ be the blowup of $C$. Then $X$ is a
Fano threefold and there is a $K$-negative extremal contraction
\cite{Mori-1982} different from $\sigma$. This contraction is
given by the birational transform of the linear system of quadrics
passing through $C$ and the fibers are birational transforms of
$2$-secant lines of $C$. This defines a $\PP^1$-bundle structure
$\varphi\colon X\to \PP^2$. We have the following
$\II$-equivariant diagram:
\[
\xymatrix{
&X\ar[dl]_{\scriptstyle{\sigma}}\ar[dr]^{\scriptstyle{\varphi}}&
\\
\PP^3\ar@{-->}[rr]&&\PP^2
}
\]
Thus $\PP^3/\II\bir X/\II$ and there is a rational curve fibration
$f\colon X/\II\to \PP^2/\II$. The action $\II$ on $\PP^2$ is
induced by an irreducible representation $\beta\colon \II /\{\pm
E\}= \AAA_5\to SL(3,\CC)$. Note that the group $\beta(\AAA_5)\cdot
\{\pm E\}$ is generated by reflections and has degrees $2,6,10$
(see \cite[Ch 17, \S 266]{Burnside-1911}). Therefore,
$\PP^2/\AAA_5$ is the weighted projective plane $\PP(2,6,10)\simeq
\PP(1,3,5)$. Let $\Delta\subset \PP^2/\AAA_5$ be the minimal curve
such that $F$ is smooth over $(\PP^2/\AAA_5)\setminus \Delta$. It
is easy to see that $\Delta$ is the image of the reflection lines
in $\PP^2$. Since there are exactly $15$ such lines (corresponding
to order $2$ elements in $\AAA_5$), the degree of $\Delta$ on
$\PP(1,3,5)$ is equal to $15$ and $\Delta$ is irreducible. Let
$\PP(1,3,5)=\operatorname{Proj} \CC[x_0,y_0,z_0]$, where $\deg x_0
= 1$, $\deg y_0 = 3$, $\deg z_0 = 5$. Using only the equality $\deg \Delta=15$
one can see that $\Delta$ is given by the following equation
\begin{multline*}
c_1z_0^3+c_2y_0^5+c_3x_0y_0^3z_0+c_4x_0^2y_0z_0^2+
c_5x_0^3y_0^4+c_6x_0^4y_0^2z_0+c_7x_0^5z_0^2+c_8x_0^6y_0^3+
\\
+c_9x_0^7y_0z_0+ c_{{10}}x_0^9y_0^2+c_{{11}}x_0
^{10}z_0+c_{{12}}x_0^{12}y_0+c_{{13}}x_0^{15} =0,
\end{multline*}
where the $c_i$ are some constants.
Consider the open set $U =\PP (1, 3, 5)\cap \{x_0\neq 0\}$. Then
in coordinates $y = y_0/x_0^3$, $z = z_0/x_0^5$ on $U\simeq
\mathbb{A}^2$ the curve $\Delta$ is defined by
\begin{multline*}
c_2y^5+c_5y^4+\left(c_3z+c_8\right)y^3+
\left(c_6z+c_{{10}}\right)y^2+ \left(c_4z^2+c_9z +c_{{12}}
\right)y+
\\
+c_1z^3+c_7z^2+c_{{11}}z+c_{{13}} =0.
\end{multline*}
Put $S_0:=U\setminus \Delta$ and $V:=f^{-1}(S_0)$. 
Then $f|_V\colon V\to S_0$ is a smooth morphism whose 
geometric fibers are isomorphic to $\PP^1$. Thus,
$V\to S_0$ is a Severi-Brauer scheme. 

\begin{lemma}
\label{lemma-Severi-Brauer}
Let $S$ be a smooth projective rational surface and let
$D\subset S$ be a reduced curve. Let $S_0:=S\setminus D$
and let $V/S_0$ be a Severi-Brauer scheme.
Assume that there is an irreducible component $D_1\subset D$
which is a smooth rational curve and such that $D_1$ meets the closure
$\overline{D-D_1}$ at a single point. Then the Severi-Brauer scheme
$V/S_0$ can be extended to $S\setminus \overline{D-D_1}$.
\end{lemma}
\begin{proof}
According to general theory, there is 1-1 correspondence between
isomorphism classes of Severi-Brauer $S_0$-schemes of relative dimension 
$n-1$ and isomorphism classes of Azumaya $\OOO_{S_0}$-algebras
of rank $n^2$. Let $A$ be the corresponding Azumaya algebra
over $S_0$. Denote by $[A]$ its class in the Brauer group 
of the function field $\Br \CC(S)$. Taking into account that 
$S$ is rational we consider the Artin-Mumford exact sequence
\cite{Artin-Mumford-1972}:
\begin{multline*}
0 \longrightarrow \Br \CC(S) \stackrel{a}{\longrightarrow} 
\bigoplus_{C\subset S} H^1(\CC(C),\QQ/\ZZ) 
\stackrel{r}{\longrightarrow}
\\
\longrightarrow \bigoplus_{P\in S} \muu^{-1} 
\stackrel{s}{\longrightarrow} \muu^{-1} \longrightarrow 0,
\end{multline*}
where $\muu^{-1}:=\cup_{n} \Hom (\muu_n,\QQ/\ZZ)$,
the first sum runs through all irreducible curves 
$C\subset S$ while the second runs through all closed points
$P\in S$ (for details we refer to \cite{Artin-Mumford-1972},
\cite{Tannenbaum-1981}). 
Note that the map $r$ ``measures'' ramification and the map $s$
is just the sum. Now let $P:=D_1\cap \overline{D-D_1}$.
Since there are no cyclic coverings of $D_1\simeq\PP^1$ ramified 
only at $P$, the $H^1(\CC(D_1),\QQ/\ZZ)$-component of 
$a([A])$ is zero, i.e., the algebra $A$ is \textit{unramified} over $D_1$. 
In this situation, $A$ can be extended over $D_1$, see, e.g., 
\cite[Prop. 6.2]{Tannenbaum-1981}. 
\end{proof}

Now we
consider the natural embedding $U = \mathbb{A}^2\hookrightarrow \PP^2
=\bar U$. Let $(x,y,z)$ be homogeneous coordinates on $\PP^2$ 
so that $\PP^2_{x,y,z}\cap \{z\neq 0\}=\mathbb{A}^2_{x,y}=U$.
Let $\bar \Delta\subset \PP^2$ be the closure of $\Delta\cap U$.
Then $\bar \Delta$ intersects the infinite line $N:=\{x=0\}$ at a single
point $P:=(0,0,1)$ which is cuspidal. By the above lemma 
the Severi-Barauer scheme $V$ can be extended to $N$. Let $L_t$ be the pencil of
lines on $\PP^2$ through $P$.
Then a general member of $L_t$ meets 
$\bar \Delta$ at $P$ and three more points.
According to \cite{Sarkisov-1982-e}
there is a \textit{standard conic bundle} $g\colon Y\to \widetilde \PP^2$
and a commutative diagram:
\[
\xymatrix{
X/\II\ar@{-->}[d]&Y\ar[d]^g\ar@{-->}[l]
\\
\PP^2&\widetilde \PP^2\ar[l]_{\psi}
}
\]
Here $\widetilde \PP^2$ is a smooth surface, 
$\psi$ is a birational morphism, and $Y\dashrightarrow X/\II$
is a birational map. By the above the discriminant curve
$\widetilde \Delta$ of $g$ is contained in $\psi^{-1}(\bar
\Delta)$. Moreover, again by Lemma \xref{lemma-Severi-Brauer}
the Severi-Barauer scheme $V$ can be extended to all exceptional 
divisors over $P$ (because $P\in \bar \Delta$ is a cuspidal point).
Thus
we may assume that $\widetilde L_t$ is a base
point free pencil such that $\widetilde L_t\cdot \widetilde
\Delta=3$. In this situation, $Y$ is rational (see
\cite{Iskovskikh-1987}), so are both $X/\II$ and $\PP^3/\II$.
\end{proof}

\subsection*{Groups of type (V)}
Since the ring of invariants
of the standard representation of $\SSS_5$
on $\CC^4$ is generated by symmetric polynomials $s_2,\dots, s_5$, 
the quotient $\PP^3/\SSS_5\simeq \PP(2,3,4,5)$ is
rational. Note also that $\AAA_5$ has a faithful three-dimensional
representation, so $\PP^3/\AAA_5$ is stably rational (more
precisely, $\PP^3/\AAA_5\times\PP^3$ is rational). T. Maeda
\cite{Maeda-T-1989} (see also \cite{Kervaire-Vust-1989}) proved
the rationality of $\mathbb A^5/\AAA_5\bir
\PP^3/\AAA_5\times\PP^2$ (over an arbitrary field). The
rationality of $\PP^3/\AAA_5$ over $\CC$ was proved in
\cite{Kervaire-Vust-1989} by an algebraic method.
Here we propose an alternative, geometric approach.

\begin{proposition}
\label{th-A5}
There is the following $\SSS_5$-equivariant diagram
\[
\xymatrix{
X\ar@{-->}[rr]^{\scriptstyle{\chi}}\ar[dr]^{\scriptstyle{\varphi_0}}
\ar[dd]_{\scriptstyle{\sigma}}
&&{X^+}\ar[ld]_{\scriptstyle{\varphi_0^+}}
\ar[dd]^{\scriptstyle{\varphi}}
\\
&{\mathcal
S_3}&
\\
\PP^3\ar@{-->}[rr]&&W
}
\]
where $\chi$ is a flop, $\varphi_0$ and $\varphi_0^+$ are small
contractions to the Segre cubic $\mathcal S_3$, $W$ is a smooth
del Pezzo surface of degree $5$ and $\varphi$ is a $\PP^1$-bundle.
\end{proposition}
As in the proof of Theorem \xref{th-SL25} the rationality of 
$\PP^3/\AAA_5$ can be proved by a
detailed analysis of the discriminant curve of the rational curve
fibration $X^+/\AAA_5\to W/\AAA_5$.

\begin{proof}
Let $G=\SSS_5$. Consider the standard representation $G=\SSS_5
\hookrightarrow GL(4,\CC)$ and the corresponding linear action on
$\PP^3$. Then $G$ permutes five points $P_1,\dots,P_5\in \PP^3$.
Let $\sigma\colon X\to \PP^3$ be the blowup of $P_1,\dots,P_5$ and
let $S_i=\sigma^{-1}(P_i)$ be the exceptional divisors. Let $H$ be
the class of hyperplane section on $\PP^3$ and let
$H^*:=\sigma^*H$. Since $\sum P_i$ is an intersection of quadric
(scheme-theoretically), the linear system $|2H^*-\sum
S_i|=|-\frac12 K_X|$ is base point free. In particular, $-K_X$ is
nef and big, i.e., $X$ is a \textit{weak Fano threefold}. It is easy to
check that $\dim |-\frac12 K_X|=4$ and the morphism
$\varphi_0\colon X\to X_0\subset \PP^4$ given by $|-\frac12 K_X|$
is birational onto its image, a three-dimensional cubic. Moreover,
$\varphi_0$ small and contracts proper transforms $L_{i,j}$ of
lines passing through $P_i$ and $P_j$. Thus $X_0$ has ten singular
points $\varphi_0(L_{i,j})$. Therefore, $X_0$ is the Segre cubic
$\mathcal S_3$ and our construction is inverse to the construction
in \S \xref{sect-Segre}, Proof of Corollary \xref{th-Segre-2}.

Since $(\Pic X)^G$ is of rank two, there is $G$-equivariant flop
$\chi\colon X \dashrightarrow X^+$. Here $X^+$ is a small
resolution of $X_0=\mathcal S_3$ obtained from $X$ by ``changing
all the signs'' \cite{Finkelnberg-1987}. There is a unique
$K$-negative $G$-extremal ray on $X^+$ \cite{Mori-1982}. Let
$\varphi\colon X^+\to W$ be its contraction.

Assume that $\rho(X^+/W)>1$. Then $\varphi$ passes through a
(non-$G$-equivariant) extremal contraction $\varphi_1\colon X^+\to
W_1$. If $\varphi_1$ is birational, then taking into account that
$-K_{X^+}$ is divisible by $2$ and the classification
\cite{Mori-1982} we get that $W_1$ is smooth and $\varphi_1$ is
the blowup of a point. Let $S_1$ be the corresponding exceptional
divisor and let $S_1,\dots,S_r$ be the $G$-orbit. By the extremal
property the $S_i$ are disjointed and give us extremal rays on
$X^+/W$. Since $\rho(X^+/W)\le 5$, $r\le 5$. On the other hand,
the image of $S_i$ on $X_0=\mathcal S_3$ is a plane. By Lemma
\xref{lemma-tran--S5} we have $r=5$ and the $S_i$ are proper
transforms of $\sigma$-exceptional divisors. Hence, 
$\chi=\mathrm{id}$ and $W=\PP^3$. Clearly, this is
impossible. Thus, $\varphi_1$ is not birational. We claim that $W$
is a surface. Indeed, $W$ cannot be a point because $-K_{X^+}$ is
not ample. Assume that $W$ is a curve and let $F$ be a general
fiber. Since $K_F=K_{X^+}|_F$ is divisible by $2$, $F\simeq
\PP^1\times \PP^1$. On the other hand the Mori cone $\Mori(X^+/W)$
has at least $5$ (non-birational) extremal rays. Each of them
gives a contraction to a surface over $W$. This is impossible
because the fibers of these contractions are contained in the fibers of
$\varphi$. Thus $W$ is a surface. By construction, the linear
system $|-\frac12 K_{X^+}|$ is the pull-back of the system
of hyperplane sections of $\mathcal S_3$. In particular, $|-\frac12
K_{X^+}|$ is base point free and a general element $D\in |-\frac12
K_{X^+}|$ is a smooth cubic surface. Since the restriction
$\varphi|_D\colon D\to W$ is birational, the surface also is
smooth and $-K_W$ is ample, i.e., $W$ is a del Pezzo surface. The
group $G$ faithfully acts on $W$. Clearly, this is not possible if
$\rho(W)\le 5$. Therefore, $\rho(W)=5$ and $K_W^2=5$. Finally,
$-K_{X^+}$ is divisible by $2$, so the fibration $\varphi$ has no
degenerate fibers. This proves the statement.
\end{proof}

\section{Monomial groups}
An action of a group $\Gamma$ on a field $\CC(x_1,\dots,x_n)$ is
said to be \textit{monomial} (with respect to $x_1,\dots, x_n$) if
for every $g\in \Gamma$ one has
\begin{equation}
\label{eq-gxi}
g(x_i)=\lambda_i(g)x_1^{m_{i,1}}\cdots x_n^{m_{i,n}},\quad
\lambda_i(g)\in \CC^*, \quad m_{i,j}\in \ZZ.
\end{equation}
Rationality of the fields of invariants of such actions was 
studied in
the series of works \cite{Haeuslein-1971}, \cite{Hajja-1983}, 
\cite{Hajja-Kang-1994}, \cite{Hajja-2000-surv}.

Any monomial action \eqref{eq-gxi} defines an integer representation
\[
\pi\colon \Gamma \to GL(n,\ZZ),\qquad g\to (m_{i,j}).
\]

Now let $V:=\CC^{n+1}$ and let $G\subset GL(V)$ be an
imprimitive group of type $(1^{n+1})$. Then $G$ permutes the
one-dimensional subspaces $V_i$ from Definition \xref{def-imp} and
contains a normal Abelian subgroup $A$ which acts diagonally in
the corresponding coordinates $x_1,\dots,x_{n+1}$. Let
$\Gamma:=G/A$. By the above, there is a natural embedding
$\Gamma\hookrightarrow \SSS_n$. Since we are assuming that the
representation $G\hookrightarrow GL(V)$ is irreducible, the group
$\Gamma\subset \SSS_{n+1}$ is transitive. Put
\[
y_1=\frac{x_1}{x_{n+1}},\dots, y_{n}=\frac{x_{n}}{x_{n+1}}\in
\CC(\PP(V)).
\]
The action of $A$ on $y_1,\dots, y_{n}$ is diagonal in these
coordinates. If $f=\sum_I a_I y_I\in \CC[y_1,\dots,y_n]$ is an
$A$-invariant, then so are all the monomials $a_Iy_I$. Hence the
ring $\CC[y_1,\dots,y_n]^A$ and its fraction field
$\CC(y_1,\dots,y_n)^A$ are generated by invariant monomials. Let
$z_1,\cdots,z_n$ be a basis of the free $\ZZ$-module
$(\text{$A$-invariant monomials in $y_i$})^*/\CC^*$. Then the
field of invariants $\CC(\PP(V))^A$ is generated by these
$z_i=z_i(y_1,\dots,y_{n})$ and they are algebraically independent.
Further, the action of $G$ on $\CC(\PP(V))^A$ is monomial with
respect to $z_1,\dots,z_{n}$. Thus the rationality question of
$\PP^n/G$ is reduced to the rationality question of the invariant
field $\CC(z_1,\dots,z_{n})^\Gamma$ of a monomial (but in
general non-linear) action.

We have two representations $\pi_{0},\, \pi\colon \Gamma \to
GL(n,\ZZ)$, corresponding monomial actions of $\Gamma$ on
$\CC(y_1,\dots,y_n)$ and $\CC(z_1,\dots,z_n)$, respectively.
Moreover, $\pi$ is a restriction of $\pi_0$ to an invariant
sublattice of finite index. In particular, $\pi\otimes
\QQ\simeq\pi_0\otimes \QQ$.

In some cases the rationality of the field
$\CC(z_1,\dots,z_{n})^\Gamma$ can be proved by purely algebraic
methods. For example, $\CC(z_1,\dots,z_{n})^\Gamma$ is rational
in the following cases:

\begin{enumerate}
\item
$\Gamma$ is cyclic of order $m$, where the class number of $m$th
cyclotomic field is $1$ \cite{Haeuslein-1971}, \cite{Hajja-1983},
\item
$\dim V\le 4$ and $G$ is meta-Abelian \cite{Hajja-1983} (in fact,
the author proved the rationality of $V/G$),
\item
$n\le 3$ and the action of $\Gamma$ on $\CC(z_1,\dots,z_n)$ is
purely monomial (i.e., all constants
$\lambda_i(g)$ in \eqref{eq-gxi} are equal to $1$)
with one exception\footnote{Professor Ming-Chang Kang
has informed me that the last unsettled case has been solved
in the paper of A. Hoshi and Y. Rikuna (to appear in ``Math. of
Computation'')} \cite{Hajja-Kang-1994}.
\end{enumerate}

From now on we consider the four-dimensional case, i.e., the case
$n=4$.
\begin{theorem}
\label{th-imprim-1-4}
Let $G\subset GL(4,\CC)$ be an imprimitive group of type $(1^4)$.
Then $\PP^3/G$ is rational.
\end{theorem}

There are the following possibilities for $\Gamma\subset \SSS_4$:
1) $\Gamma$ is a cyclic group of order $4$, 2) Klein group
$\mathfrak V_4$, 3) $\Gamma$ is a dihedral group $\DDD_4$ of order
$8$, 4) $\Gamma=\AAA_4$, and 5) $\Gamma=\SSS_4$. By the above it
is sufficient to prove the rationality of the field
$\CC(z_1,z_2,z_3)^\Gamma$. The theorem was proved in the unpublished
manuscript of I. Kolpakov-Miroshnichenko and the author (1986) by
case by case consideration of the action of $\Gamma$ on
$\CC(z_1,z_2,z_3)$. We consider here only the case
$\Gamma=\AAA_4$. Case $\Gamma=\SSS_4$ can be treated in a similar
way. Cases $\Gamma=\ZZ_4$, $\mathfrak V_4$, and $\DDD_4$ are
easier. Moreover, in the these three cases $G$ is also imprimitive
of type $(2^2)$ and then the rationality can be proved also by
another method, see \S \xref{sect-imp-2-2}.

The group $\AAA_4$ is generated by two elements: $\delta=(1,2)(3,4)$ and
$\theta=(1,2,3)$. One has:
\[
\pi_0(\theta)= \left(
\begin{array}{rrr}
0&1&0
\\
0&0&1
\\
1&0&0
\end{array} \right)
\qquad
 \pi_0(\delta)= \left(
\begin{array}{rrr}
0&1&-1
\\
1&0&-1
\\
0&0&-1
\end{array} \right)
\]
Using the classification of finite subgroups in
$GL(3,\ZZ)$ \cite{Tahara} we get the following possibilities:
\begin{itemize}
\item[$\Gamma_9^{12}$:]
\quad $\pi(\theta)= \left(
\begin {array}{rrr} 0&1&0
\\
0&0&1
\\
1&0&0
\end {array}\right)
 \quad \pi(\delta)= \left(
\begin {array}{rrr} 
-1&\phantom{-}0&\phantom{-}0
\\0&1&0
\\0&0&-1
\end {array}\right)$

\item[$\Gamma_{10}^{12}$:]
\quad $\pi(\theta)= \left(
\begin {array}{rrr}
0&1&0
\\
0&0&1
\\
1&0&0
\end {array}\right)
 \quad \pi(\delta)=
\left(\begin {array}{rrr} 0&-1&1
\\
\phantom{-}0&-1&\phantom{-}0
\\
1&-1&0
\end {array}\right)$

\item[$\Gamma_{11}^{12}$:]
\quad $\pi(\theta)= \left(\begin {array}{rrr} 0&1&0
\\
0&0&1
\\
1&0&0
\end {array}\right)\quad
\pi(\delta)= \left(\begin {array}{rrr} -1&-1&-1\\0&0&1
\\0&1&0
\end {array}\right)$
\end{itemize}

\subsection*{Case $\Gamma_9^{12}$} 
By \eqref{eq-gxi} the action of $\Gamma=\AAA_4$ on
$\CC(z_1,z_2,z_3)$ has the following form
\[
(a_1 z_2,a_2 z_3,a_3 z_1) \stackrel{\theta}{\longleftarrow}
(z_1,z_2,z_3) \stackrel{\delta}{\longrightarrow} (b /z_1,c
z_2,b'/z_3),\quad a_i, b, b',c \in \CC^*.
\]
After coordinate change of the form $z_i\to \lambda_i z_i$,
$\lambda_i\in \CC^*$ we may assume that $a_1=a_2=a_3=b'=1$. Since
$\delta^2=1$, $c=\pm 1$. From other relations between $\theta$ and
$\delta$ we get $b=c=\pm 1$. Therefore,
\[
(z_2,z_3,z_1) \stackrel{\theta}{\longleftarrow}
(z_1,z_2,z_3)\stackrel{\delta}{\longrightarrow} (b/z_1,b
z_2,1/z_3), \qquad b=\pm 1.
\]
Assume that $b=c=1$. Then after the coordinate change of the form
\[
z_1'=(z_1+1)/(z_1-1),\qquad z_2'=(z_2+1)/(z_2-1),\qquad
z_3'=(z_3+1)/(z_3-1)
\]
the action will be linear:
\[
(z_2',z_3',z_1') \stackrel{\theta}{\longleftarrow}
(z_1',z_2',z_3')\stackrel{\delta}{\longrightarrow}
(-z_1',z_2',-z_3')
\]
By Proposition \xref{prop-first}  the field
$\CC(z_1',z_2',z_3')^\Gamma$ is rational.

Now assume that $b=c=-1$. Regard $z_1$, $z_2$, $z_3$ as
non-homogeneous coordinates on $\PP^1\times\PP^1\times\PP^1$. We
get an action of $\Gamma=\AAA_4$ on $\PP^1\times\PP^1\times\PP^1$
by regular automorphisms. Consider the Segre embedding
$\PP^1\times\PP^1\times\PP^1 \hookrightarrow \PP^7$
\[
(z_1, z_2, z_3) \longrightarrow (z_1z_2z_3,\ z_1,\
z_2,\ z_{3},\ z_2z_3,\ z_1z _3,\
z_1z_2,\ 1)
\]
and let $t_1,\dots,t_8$ be the corresponding coordinates in
$\PP^7$. This induces the following representation of the
tetrahedron group $\TT$ into $GL(8,\CC)$:
\begin{multline*}
(t_1,\ t_3,\ t_4,\ t_2,\ t_6,\ t_7,\ t_5,\ t_8)
\stackrel{\widetilde \theta}{\longleftarrow} (t_1,\dots,t_8)
\stackrel{\widetilde \delta}{\longrightarrow} 
\\
(t_3,\ -t_4,\ -t_1,\
t_2,\ -t_7,\ -t_8,\ t_5,\ t_6).
\end{multline*}
It is easy to check that this representation is the direct sum of
four two-dimensional faithful representations. Therefore, there
are four $\Gamma$-invariant lines $L_i$ in $\PP^7$. Clearly
$\PP^1\times\PP^1\times\PP^1$ contains no $\AAA_4$-invariant
lines. On the other hand, the intersection
$\PP^1\times\PP^1\times\PP^1\cap L_i$ consists of at most two
points (because the image of the Segre embedding is an
intersection of quadrics). This immediately implies that all the
$L_i$ are disjointed from $\PP^1\times\PP^1\times\PP^1$. Let
$V\simeq \PP^5$ be the linear span of $L_2$, $L_3$, $L_4$ and let
$\mathcal H$ be the pencil of hyperplane sections of
$\PP^1\times\PP^1\times\PP^1$ passing through $V$.

We claim that the general member of $H\in \mathcal H$ is smooth.
Let $B=V\cap\PP^1\times\PP^1\times\PP^1$ be the base locus of
$\mathcal H$. By Bertini's theorem $\Sing (H)\subset B$. 
If $H$ is not normal, then by the adjunction formula $H$ 
is singular along a line. 
On the other hand, $\PP^1\times\PP^1\times\PP^1$ does not contain any
$\AAA_4$-invariant lines. This immediately implies that $H$ is normal. Since $-K_H$
is ample, $H$ is either a cone over an elliptic curve or a del
Pezzo surface with Du Val singularities. In both cases the number
of singular points is at most two. On the other hand, $\Sing(H)$
coincides with the set of points where the dimension of
Zariski
tangent space $T_{P,B}$ jumps. Therefore, $\Sing(H)$ is
$\AAA_4$-invariant. This immediately implies the existence of
$\AAA_4$-invariant point $P\in \PP^1\times\PP^1\times\PP^1$, a
contradiction.

Thus, a general member $H\in \mathcal H$ is a smooth del Pezzo surface
of degree $6$. Now consider the projection of
$\PP^1\times\PP^1\times\PP^1$ from $V$ to $L_1$. By blowing up $B$
we obtain a $\AAA_4$-equivariant fibration $f\colon Y\to \PP^1$
whose general fiber is a smooth del Pezzo surface of degree $6$.
This induces a fibration $Y/\AAA_4\to \PP^1/\AAA_4$ with the same
type of general fiber. Such a fibration is rational over
$\CC(\PP^1/\AAA_4)$ (see \cite[Ch. 4]{Manin-Cubic-forms-e-II}), so
$(\PP^1\times\PP^1\times\PP^1)/\AAA_4\bir Y/\AAA_4$ is also
rational.

\subsection*{Case $\Gamma_{10}^{12}$}
Regard $z_1$, $z_2$, $z_3$ as non-homogeneous coordinates on
$\PP^3$. We get a linear action of $\Gamma=\AAA_4$ on $\PP^3$.
This action is either reducible or imprimitive. In the first case,
$\CC(\PP(V))^G$ is rational. In the second one, we can repeat
inductive procedure to reduce the problem to the smaller order of
$\Gamma$.

\subsection*{Case $\Gamma_{11}^{12}$}
Then after coordinate change the action of $\Gamma=\AAA_4$ on
$\CC(z_1,z_2,z_3)$ has the following form
\[
(z_2, z_3, z_1) \stackrel{\theta}{\longleftarrow} (z_1,z_2,z_3)
\stackrel{\delta}{\longrightarrow} (1/(z_1z_2z_3), z_3,z_2).
\]
i.e., it is purely monomial. The field of invariants is rational
by \cite{Hajja-Kang-1994}.

\section{Imprimitive case $(2^2)$}
\label{sect-imp-2-2}
\begin{theorem}
Let $G\subset GL(4,\CC)$ be an imprimitive group of type $(2^2)$.
Then $\PP^3/G$ is rational.
\end{theorem}
\begin{proof}[Outline of the proof]
Let $G\subset SL(V)$, $\dim V=4$ be an imprimitive group of type
$(2^2)$. Then there is a decomposition $V=V_1\oplus V_2$, $\dim
V_i=2$ and a subgroup $N\subset G$ of index $2$ such that $N\cdot
V_i=V_i$. We may assume that $N$ is not Abelian (otherwise $G$ is
reducible). 

The decomposition $V=V_1\oplus V_2$ defines two skew lines $L_1$,
$L_2$ in $\PP^3=\PP(V)$. Let $H$ be the class of hyperplane in
$\PP^3$. The image of the map $\PP^3\dashrightarrow \PP^3$ given
by the linear system $2H-L_1-L_2$ is the two-dimensional quadric
$\PP^1\times \PP^1$ and fibers of this map are lines meeting both
$L_1$ and $L_2$. We have the following $G$-equivariant diagram:
\[
\xymatrix{
&{\widetilde \PP^3}\ar[ld]_{\scriptstyle{\sigma}}\ar[rd]^{\scriptstyle{\varphi}}&
\\
\PP^3\ar@{-->}[rr]&&{\PP^1\times\PP^1}
}
\]
where $\sigma\colon \widetilde \PP^3\to \PP^3$ is the blow up of
$L_1\cup L_2$ and $\varphi$ is a $\PP^1$-bundle. Let $S_i$ be the
corresponding exceptional divisors. There is a rational curve
fibration $f\colon \widetilde \PP^3/ G\to S=(\PP^1\times \PP^1)/
G$. The rationality of $\PP^3/ G$ follows from detailed analysis
of the action of $G$ on $\PP^1\times \PP^1$ (cf. Proof of Theorem
\xref{th-SL25}).

Consider for example the case when restrictions $N\to GL(V_i)$
are injective. Let $\pi\colon \PP^1\times\PP^1\to (\PP^1\times
\PP^1)/ G$ be the quotient map and let $\Delta\subset S$ be the
discriminant curve. Then, $B:=\pi^{-1}(\Delta)$ is contained in the
ramification divisor, the union of one-dimensional components of
the locus of points with non-trivial stabilizer. Take a point
$P=(x,y)\in B$. If $g\in \St(P)$, then $g(x)=x$, $g(y)=y$. 
If $g\in N$ or
$g^2\neq \lambda E$, there
are only a finite number of such points. 
This implies that the ramification index over
each component of
$\Delta$ is equal to $2$. Consider the pencil $C_t$ of $(1,0)$-curves 
of $\PP^1\times \PP^1$.
Let $D_t:=\pi(C_t)$. Then $D_t$ is a base point free pencil. Hence,
\[
K_S\cdot D_t=2p_a(D_t)-2-D_t^2=-2. 
\]
Further, by the Hurwitz formula,
\[
K_{\PP^1\times\PP^1}=\pi^*\left(K_S+\frac 12 \Delta\right).
\]
This yields
\begin{multline*}
0>\frac 2{\deg \pi} K_{\PP^1\times\PP^1}\cdot \pi^* D_t=
\frac 2{\deg \pi}\pi^*\left(K_S+\frac 12 \Delta\right)\cdot \pi^* D_t=
\\
=2\left(K_S+\frac 12 \Delta\right)\cdot D_t =\Delta \cdot D_t-4.
\end{multline*}
Therefore, $\Delta\cdot D_t\le 3$. Finally, as in the proof of
Theorem \xref{th-SL25} we deduce the rationality of $\widetilde
\PP^3/G$ using \cite{Sarkisov-1982-e} and \cite{Iskovskikh-1987}.
\end{proof}

\section{Final remarks and open questions}
As a consequence of the above results we have the following
\begin{theorem}
Let $G\subset GL(4,\CC)$ a finite subgroup. Assume that $G$ is
solvable. Then $\PP^3/G$ is rational.
\end{theorem}

\subsection*{Remaining cases}
The rationality question for $\PP^3/G$ and $\CC^4/G$ is still open
for the following groups (up to scalar multiplication).
For convenience of the reader we give a short description of 
the group action. 

\subsubsection*{Type (I)} 
The rationality of $\PP^3/G$ is unknown only for 
$G=\Psi(\OO,\II)$.
Note that $\PP^3/\Psi(\OO,\II)$ is stably rational,
see Theorem \xref{th-O-I-st-rat}.

\subsubsection*{Types (VI), (VII) and (VIII)}
Unsolved cases are 
$\widetilde \SSS_5$, $\widetilde \AAA_6$, $\widetilde \SSS_6$, and
$\widetilde \AAA_7$.
The corresponding actions on $\PP^3$ are given by projective representations
of $\SSS_5$, $\AAA_6$, $\SSS_6$, and
$\AAA_7$ into $PGL(4,\CC)$. Recall that the Schur multiplier $H^2(\SSS_n,\CC^*)$ 
of the symmetric group is isomorphic to 
$\muu_2$ for $n\ge 4$. Therefore, any 
projective representation of $\SSS_n$, $n\ge 4$ is induced by a
linear representation of a central extension $\tilde \SSS_n$
by $\muu_2$.

Following Schur 
\cite{Schur-1911}, \cite{Schur-1911-2001} we 
give an explicit matrix representation $\tilde \SSS_6\to
GL(4,\CC)$. 
Consider the following matrices
\[
E:=\begin{pmatrix}
1&0
\\
0&1
\end{pmatrix}
\quad
A:=\begin{pmatrix}
0&\sqrt{-1}
\\
\sqrt{-1}&0
\end{pmatrix}
\quad 
B:= 
\begin{pmatrix}
0&-1
\\
1&0
\end{pmatrix}
\quad 
C:=
\begin{pmatrix}
1&0
\\
0&-1
\end{pmatrix}.
\]
Next take the following Kronecker products
\[
M_1:=C\otimes A,\ M_2:=C\otimes B, \ 
M_3:=A\otimes E,\ M_4:=B\otimes E, \
M_5:=\sqrt{-1}C^{\otimes 2}.
\]
One can easily check the relations 
\begin{equation}
\label{eq-Spin}
M_j^2=-E_4, \quad M_jM_k=-M_kM_j,\qquad 1\le j\neq k\le 5, 
\end{equation}
where $E_4$ is the identity $4\times 4$ matrix.
Now put
\begin{equation}
\label{eq-Spin-Tk-def}
T_k:=\frac{1}{\sqrt{2k}}\left( -\sqrt{k-1}M_{k-1}+\sqrt{k+1}M_{k}\right),
\quad k=1,\dots,5.
\end{equation}
From \eqref{eq-Spin} we have 
\begin{equation}
\label{eq-Spin-Tk}
\begin{array}{l}
T_k^2=-E_4,\qquad (T_kT_{k+1})^3=-E_4,
\\[10pt]
T_jT_k=-T_kT_j\quad \text{for}\quad k>j+1.
\end{array}
\end{equation}
Generators $T_j$ and relations \eqref{eq-Spin-Tk}
determine an abstract group $\tilde \SSS_6$
that is the central extension of $\SSS_6$ by $\muu_2$.
Here $T_j$ corresponds to the transposition interchanging $j$ 
and $j+1$.
In fact, the constructed representation 
$\tilde \SSS_6\to GL(4,\CC)$ is obtained from 
the standard action of $\SSS_6$ on the Clifford algebra 
$A(\CC^5)$ (cf., e.g., \cite{Chevalley-1954-a}). 
Taking compositions with embeddings into $\SSS_6$
we get also projective representations 
$\SSS_5\hookrightarrow PGL(4,\CC)$ and $\AAA_6\hookrightarrow PGL(4,\CC)$.

For $\tilde \AAA_7$, similarly put 
\begin{multline*}
M_1:=C^{\otimes 2}\otimes A,\ M_3:=C\otimes A\otimes E,\ 
M_5:=A\otimes E^{\otimes 2},\ 
\\
M_2:=C^{\otimes 2}\otimes B,\ M_4:=C\otimes B\otimes E, \ 
M_6:=B\otimes E^{\otimes 2},\
M_7:=\sqrt{-1}C^{\otimes 3}
\end{multline*}
and define $T_k$ by formula \eqref{eq-Spin-Tk-def}.
As above, we get an irreducible linear representation $\tilde S_7\to GL(8,\CC)$.
The restriction to $\tilde \AAA_7$ splits as a direct sum of 
two $4$-dimensional faithful representations. 
The explicit generators of $\AAA_7\subset PGL(4,\CC)$ 
can be taken as follows 
\[
S= \begin{pmatrix}
1&0&0&0
\\
0&\beta&0&0
\\
0&0&{\beta}^{4}&0
\\
0&0&0&{\beta} ^{2}
\end{pmatrix}
\quad W= \frac{1}{\sqrt {-7}}
\begin{pmatrix}
{p}^{2}&1&1&1
\\
1&-p&-q&-p
\\
1&-p&-p&-q
\\
1&-q&-p&-p
\end{pmatrix}
\]
where $\beta$ is a primitive $7$-th root of unity,
$p:=\beta+\beta^2+\beta^4$, and $q:=\beta^3+\beta^5+\beta^6$, 
see \cite{Blichfeldt}. 
Then the isomorphism between the subgroup in $PGL(4,\CC)$ generated by 
$S$, $W$ and $\AAA_7$ is given by 
\[
S\to (1,2,3,4,5,6,7),\qquad W\to (2,3,5)(4,6,7).
\]

\begin{remark}
The group $\AAA_6$ has also a central extension $\hat \AAA_6$
by $\muu_3$ which admits a $3$-dimensional representation 
$\delta\colon \hat \AAA_6\to SL(3,\CC)$. The group $\delta(\hat \AAA_6)$
is projectively equivalent to the complex reflection group 
of order $360\cdot 6$ \cite{Shephard-Todd}, \cite{Cohen-1976}. 
Possibly this can be used to prove the stable rationality
of $\PP^3/\tilde \AAA_6$. 
\end{remark}

\subsubsection*{Type (IX), $G\simeq SL(2,\FF_7)$}
The representation $G\hookrightarrow SL(4,\CC)$ can be described as follows.
Fix a character $\chi \colon \FF_7^*\to \CC^*$ 
and consider the following $\CC$-vector space of 
$\CC$-valued functions on 
$\FF_7^2\setminus \{0\}$
\[
V_{\chi}:=\left\{f\colon \left(\FF_7^2\setminus \{0\}\right)\to \CC \mid 
f(\lambda x)= \chi(\lambda)f(x),\ \forall \lambda\in \FF_7^*\right\}.
\]
It is easy to see that $\dim_\CC V_\chi=8$ and $SL(2,\FF_7)$ naturally acts on $V_\chi$.
Thus we have a representation $\rho_\chi\colon SL(2,\FF_7)\to SL(V_\chi)=SL(8,\CC)$.
If $\chi^2=1$, $\chi\neq 1$, then $\rho_\chi$ splits as a direct sum of 
two $4$-dimensional faithful representations. This gives as a subgroup 
in $SL(4,\CC)$ isomorphic to $SL(2,\FF_7)$. Explicit matrices 
can be found, e.g., in \cite{Blichfeldt}.
The ring of invariants is completely described, see 
\cite{Mallows-Sloane-1973}.
Note that $SL(2,\FF_7)$ has a $3$-dimensional non-faithful representation
$\delta\colon SL(2,\FF_7)\to SL(3,\CC)$. 
The group $G':=G\times \{\pm E\}$ is generated by complex
reflections and the degrees are $4$, $6$ and $14$
(see \cite[\S 139]{Weber-1896}, \cite{Cohen-1976}).

\subsubsection*{Type (XI)}
The rationality of $\PP^3/G$ is unknown only for 
two groups of order $64\cdot 360$ and $64\cdot 60$, see
Corollaries \xref{th-Segre-2} and \xref{th-Segre-3}.
Using Theorem \xref{th-Segre-1} one can easily get equations for birational models of 
these quotients in terms of discriminants.

\subsection*{$p$-groups}
It is known that the answer to Noether's problem is negative in
higher dimensions: there are examples of $p$-groups such that
$\Bbbk(\{x_g\}_{g\in G})^G$ is not rational (and even not stably
rational) \cite{Saltman-1984}, \cite{Shafarevich-1990-re}.
Moreover, for any prime $p$ there is a group $G$ of order $p^6$
such that $\Bbbk(\{x_g\}_{g\in G})^G$ is not stably rational
\cite{Bogomolov-1987-e}. On the other hand, it is known that the
field of invariant of any linear action of a $p$-group of order
$\le p^4$ on $\Bbbk(x_1,\dots,x_N)$ is rational whenever $\Bbbk$
contains a primitive $p^e$-th root of unity, where $p^e$ is the
exponent of $G$ \cite{Chu-Kang-2001}, see also
\cite{Beneish-2003}. The following question is open:
\begin{quote}
is it true that for any group $G$ of order $p^5$ and any linear action 
of $G$ on
$\CC^n$ the quotient $\CC^n/G$ is stably rational?
\end{quote}
It seems that the answer is positive (the author checked this for $p=2$).
Note that in characteristic $p>0$ any linear action of $p$-group
is rational \cite{Kuniyoshi-1955}.

 \def\cprime{$'$} \def\polhk#1{\setbox0=\hbox{#1}{\ooalign{\hidewidth
  \lower1.5ex\hbox{`}\hidewidth\crcr\unhbox0}}}

%\bibliography{/home/prokhoro/bib/my_ref,/home/prokhoro/bib/specific/inv}

\begin{thebibliography}{BCTSSD85}

\bibitem[AM72]{Artin-Mumford-1972}
M.~Artin and D.~Mumford.
\newblock Some elementary examples of unirational varieties which are not
  rational.
\newblock {\em Proc. London Math. Soc. (3)}, 25:75--95, 1972.

\bibitem[BCTSSD85]{Beauville-Colliot-Thelene-Sansuc-1985}
Arnaud Beauville, Jean-Louis Colliot-Th{\'e}l{\`e}ne, Jean-Jacques Sansuc, and
  Peter Swinnerton-Dyer.
\newblock Vari\'et\'es stablement rationnelles non rationnelles.
\newblock {\em Ann. of Math. (2)}, 121(2):283--318, 1985.

\bibitem[Ben03]{Beneish-2003}
Esther Beneish.
\newblock Stable rationality of certain invariant fields.
\newblock {\em J. Algebra}, 269(2):373--380, 2003.

\bibitem[Bli17]{Blichfeldt}
H.~F. Blichfeldt.
\newblock {\em Finite collineation groups}.
\newblock The Univ. Chicago Press, Chicago, 1917.

\bibitem[Bog88]{Bogomolov-1987-e}
F.~A. Bogomolov.
\newblock The {B}rauer group of quotient spaces of linear representations.
\newblock {\em Math. USSR-Izv.}, 30(3):455--485, 1988.

\bibitem[Bur55]{Burnside-1911}
W.~Burnside.
\newblock {\em Theory of groups of finite order}.
\newblock Dover Publications Inc., New York, 1955.
\newblock 2d ed.

\bibitem[Che54]{Chevalley-1954-a}
Claude~C. Chevalley.
\newblock {\em The algebraic theory of spinors}.
\newblock Columbia University Press, New York, 1954.

\bibitem[Che55]{Chevalley-1955}
C.~Chevalley.
\newblock Invariants of finite groups generated by reflections.
\newblock {\em Amer. J. Math}, 77:778--782, 1955.

\bibitem[CK01]{Chu-Kang-2001}
Huah Chu and Ming-chang Kang.
\newblock Rationality of {$p$}-group actions.
\newblock {\em J. Algebra}, 237(2):673--690, 2001.

\bibitem[Coh76]{Cohen-1976}
Arjeh~M. Cohen.
\newblock Finite complex reflection groups.
\newblock {\em Ann. Sci. \'Ecole Norm. Sup. (4)}, 9(3):379--436, 1976.

\bibitem[CTS05]{colliotthelene-Sansuc-2005}
Jean-Louis Colliot-Th\'el\`ene and Jean-Jacques Sansuc.
\newblock The rationality problem for fields of invariants under linear
  algebraic groups (with special regards to the {B}rauer group).
\newblock Preprint math/0507154, 2005.

\bibitem[DO88]{Dolgachev-Ortland-1989}
Igor Dolgachev and David Ortland.
\newblock Point sets in projective spaces and theta functions.
\newblock {\em Ast\'erisque}, (165):210 pp. (1989), 1988.

\bibitem[Dol87]{Dolgachev-1987}
Igor~V. Dolgachev.
\newblock Rationality of fields of invariants.
\newblock In {\em Algebraic geometry, Bowdoin, 1985 (Brunswick, Maine, 1985)},
  volume~46 of {\em Proc. Sympos. Pure Math.}, pages 3--16. Amer. Math. Soc.,
  Providence, RI, 1987.

\bibitem[EM73]{Endo-Miyata-1973}
Shizuo End{\^o} and Takehiko Miyata.
\newblock Invariants of finite abelian groups.
\newblock {\em J. Math. Soc. Japan}, 25:7--26, 1973.

\bibitem[Fei71]{Feit-1971}
Walter Feit.
\newblock The current situation in the theory of finite simple groups.
\newblock In {\em Actes du Congr\`es International des Math\'ematiciens (Nice,
  1970), Tome 1}, pages 55--93. Gauthier-Villars, Paris, 1971.

\bibitem[Fin87]{Finkelnberg-1987}
Hans Finkelnberg.
\newblock Small resolutions of the {S}egre cubic.
\newblock {\em Nederl. Akad. Wetensch. Indag. Math.}, 49(3):261--277, 1987.

\bibitem[Fis15]{Fischer-E-1915-a}
E.~Fischer.
\newblock Die {I}somorphie der {I}nvariantenkorper der endlichen {A}belschen
  {G}ruppen linearer {T}ransformationen.
\newblock {\em Nachr. Konigl. Ges. Wiss. Gottingen}, pages 77--80, 1915.

\bibitem[FN89]{Furushima-Nakayama-1989}
Mikio Furushima and Noboru Nakayama.
\newblock A new construction of a compactification of {${\bf C}\sp 3$}.
\newblock {\em Tohoku Math. J. (2)}, 41(4):543--560, 1989.

\bibitem[For84]{Formanek-1984}
E.~Formanek.
\newblock Rational function fields. {N}oether's problem and related questions.
\newblock {\em J. Pure Appl. Algebra}, 31(1-3):28--36, 1984.

\bibitem[Hae71]{Haeuslein-1971}
G.~K. Haeuslein.
\newblock On the invariants of finite groups having an abelian normal subgroup
  of prime index.
\newblock {\em J. London Math. Soc. (2)}, 3:355--360, 1971.

\bibitem[Haj83]{Hajja-1983}
Mowaffaq Hajja.
\newblock Rational invariants of meta-abelian groups of linear automorphisms.
\newblock {\em J. Algebra}, 80(2):295--305, 1983.

\bibitem[Haj00]{Hajja-2000-surv}
Mowaffaq Hajja.
\newblock Linear and monomial automorphisms of fields of rational functions:
  some elementary issues.
\newblock In {\em Algebra and number theory (Fez)}, volume 208 of {\em Lecture
  Notes in Pure and Appl. Math.}, pages 137--148. Dekker, New York, 2000.

\bibitem[HK94]{Hajja-Kang-1994}
Mowaffaq Hajja and Ming~Chang Kang.
\newblock Three-dimensional purely monomial group actions.
\newblock {\em J. Algebra}, 170(3):805--860, 1994.

\bibitem[IP99]{Iskovskikh-Prokhorov-1999}
V.~A. Iskovskikh and Yu.~G. Prokhorov.
\newblock {\em Fano varieties. {A}lgebraic geometry. {V}.}, volume~47 of {\em
  Encyclopaedia Math. Sci.}
\newblock Springer, Berlin, 1999.

\bibitem[Isk87]{Iskovskikh-1987}
V.~A. Iskovskikh.
\newblock On the rationality problem for conic bundles.
\newblock {\em Duke Math. J.}, 54:271--294, 1987.

\bibitem[Jor70]{Jordan-1870}
C.~Jordan.
\newblock {\em Trait{\'e} des substitutions des {\'e}quations alg{\'e}briques}.
\newblock Paris, 1870.

\bibitem[KM99]{Kemper-1999}
Gregor Kemper and Gunter Malle.
\newblock Invariant fields of finite irreducible reflection groups.
\newblock {\em Math. Ann.}, 315(4):569--586, 1999.

\bibitem[KMP87]{Kolpakov-Prokhorov-1987}
I.~Ya. Kolpakov-Miroshnichenko and Yu.~G. Prokhorov.
\newblock Rationality of the field of invariants of the faithful
  four-dimensional representation of the icosahedral group.
\newblock {\em Mat. Zametki}, 41(4):479--483, 619, 1987.

\bibitem[KMP92]{Kolpakov-Prokhorov-1992}
I.~Ya. Kolpakov-Miroshnichenko and Yu.~G. Prokhorov.
\newblock Rationality construction of fields of invariants of some finite
  four-dimensional linear groups associated with {F}ano threefolds.
\newblock {\em Math. Notes}, 51:74--76, 1992.

\bibitem[KMP93]{Kolpakov-Prokhorov-1993}
I.~Ya. Kolpakov-Miroshnichenko and Yu.~G. Prokhorov.
\newblock Rationality of fields of invariants of some four-dimensional linear
  groups{,} and an equivariant construction related to the {S}egre cubic.
\newblock {\em Math. USSR Sbornik}, 74(1):169--183, 1993.

\bibitem[Koi03]{Koike}
Kenji Koike.
\newblock Remarks on the {S}egre cubic.
\newblock {\em Arch. Math. (Basel)}, 81(2):155--160, 2003.

\bibitem[Kun55]{Kuniyoshi-1955}
Hideo Kuniyoshi.
\newblock On a problem of {C}hevalley.
\newblock {\em Nagoya Math. J.}, 8:65--67, 1955.

\bibitem[KV89]{Kervaire-Vust-1989}
Michel Kervaire and Thierry Vust.
\newblock Fractions rationnelles invariantes par un groupe fini: quelques
  exemples.
\newblock In {\em Algebraische Transformationsgruppen und Invariantentheorie},
  volume~13 of {\em DMV Sem.}, pages 157--179. Birkh\"auser, Basel, 1989.

\bibitem[Len74]{Lenstra-1974}
H.~W. Lenstra, Jr.
\newblock Rational functions invariant under a finite abelian group.
\newblock {\em Invent. Math.}, 25:299--325, 1974.

\bibitem[Mae89]{Maeda-T-1989}
T.~Maeda.
\newblock Noether's problem for {$A\sb 5$}.
\newblock {\em J. Algebra}, 125(2):418--430, 1989.

\bibitem[Man86]{Manin-Cubic-forms-e-II}
Yu.~I. Manin.
\newblock {\em Cubic forms}, volume~4 of {\em North-Holland Mathematical
  Library}.
\newblock North-Holland Publishing Co., Amsterdam, second edition, 1986.
\newblock Algebra, geometry, arithmetic, Translated from the Russian by M.
  Hazewinkel.

\bibitem[Mil58]{Miller-1958}
Donald~W. Miller.
\newblock On a theorem of {H}\"older.
\newblock {\em Amer. Math. Monthly}, 65:252--254, 1958.

\bibitem[Miy71]{Miyata-K-1971}
Takehiko Miyata.
\newblock Invariants of certain groups. {I}.
\newblock {\em Nagoya Math. J.}, 41:69--73, 1971.

\bibitem[Mor82]{Mori-1982}
S.~Mori.
\newblock Threefolds whose canonical bundles are not numerically effective.
\newblock {\em Ann. Math.}, 115:133--176, 1982.

\bibitem[MS73]{Mallows-Sloane-1973}
C.~L. Mallows and N.~J.~A. Sloane.
\newblock On the invariants of a linear group of order $336$.
\newblock {\em Proc. Camb. Phil. Soc.}, 74:435--440, 1973.

\bibitem[MU83]{Mukai-Umemura-1983}
Shigeru Mukai and Hiroshi Umemura.
\newblock Minimal rational threefolds.
\newblock In {\em Algebraic geometry (Tokyo/Kyoto, 1982)}, volume 1016 of {\em
  Lecture Notes in Math.}, pages 490--518. Springer, Berlin, 1983.

\bibitem[Nak87]{Nakamura-1987}
Iku Nakamura.
\newblock Moishezon threefolds homeomorphic to {${\bf P}\sp 3$}.
\newblock {\em J. Math. Soc. Japan}, 39(3):521--535, 1987.

\bibitem[Nak89]{Nakano-1989b}
Tetsuo Nakano.
\newblock On equivariant completions of {$3$}-dimensional homogeneous spaces of
  {${\rm SL}(2,{\bf C})$}.
\newblock {\em Japan. J. Math. (N.S.)}, 15(2):221--273, 1989.

\bibitem[Noe13]{Noether_1913}
E.~Noether.
\newblock Rationale {F}unctionenk{\"o}rper.
\newblock {\em Jahbericht Deutsh. Math. Verein.}, 22:316--319, 1913.

\bibitem[Sal84]{Saltman-1984}
David~J. Saltman.
\newblock Noether's problem over an algebraically closed field.
\newblock {\em Invent. Math.}, 77(1):71--84, 1984.

\bibitem[Sar82]{Sarkisov-1982-e}
V.~G. Sarkisov.
\newblock On the structure of conic bundles.
\newblock {\em Math. USSR, Izv.}, 120:355--390, 1982.

\bibitem[Sch11]{Schur-1911}
I.~Schur.
\newblock {\"U}ber die {D}arstellung der symmetrischen und der alternierenden
  {G}ruppe durch gebrochene lineare {S}ubstitutionen.
\newblock {\em J. Reine Angew. Math.}, 139:155--250, 1911.

\bibitem[Sch01]{Schur-1911-2001}
J.~Schur.
\newblock On the representation of the symmetric and alternating groups by
  fractional linear substitutions.
\newblock {\em Internat. J. Theoret. Phys.}, 40(1):413--458, 2001.
\newblock Translated from the German [J.\ Reine Angew.\ Math.\ {\bf 139}
  (1911), 155--250] by Marc-Felix Otto.

\bibitem[Seg63]{Segre-1963}
C.~Segre.
\newblock {\em Opere}, volume~4.
\newblock Edizioni Cremoneze, Roma, 1963.

\bibitem[Sha90]{Shafarevich-1990-re}
I.~R. Shafarevich.
\newblock The {L}\"uroth problem.
\newblock {\em Trudy Mat. Inst. Steklov.}, 183:199--204, 229, 1990.
\newblock Translated in Proc.\ Steklov Inst.\ Math.\ {\bf 1991}, no.\ 4,
  241--246, Galois theory, rings, algebraic groups and their applications
  (Russian).

\bibitem[Spr77]{Springer-1977}
T.~A. Springer.
\newblock {\em Invariant theory}.
\newblock Springer-Verlag, Berlin, 1977.
\newblock Lecture Notes in Mathematics, Vol. 585.

\bibitem[SR49]{Semple-Roth-1949}
J.~G. Semple and L.~Roth.
\newblock {\em Introduction to {A}lgebraic {G}eometry}.
\newblock Oxford, at the Clarendon Press, 1949.

\bibitem[ST54]{Shephard-Todd}
G.~C. Shephard and J.~A. Todd.
\newblock Finite unitary reflection groups.
\newblock {\em Canadian J. Math.}, 6:274--304, 1954.

\bibitem[Swa69]{Swan-1969}
Richard~G. Swan.
\newblock Invariant rational functions and a problem of {S}teenrod.
\newblock {\em Invent. Math.}, 7:148--158, 1969.

\bibitem[Tah71]{Tahara}
Ken~ichi Tahara.
\newblock On the finite subgroups of {${\rm GL}(3,\,Z)$}.
\newblock {\em Nagoya Math. J.}, 41:169--209, 1971.

\bibitem[Tak89]{Takeuchi-1989}
Kiyohiko Takeuchi.
\newblock Some birational maps of {F}ano {$3$}-folds.
\newblock {\em Compositio Math.}, 71(3):265--283, 1989.

\bibitem[Tan81]{Tannenbaum-1981}
Allen Tannenbaum.
\newblock The {B}rauer group and unirationality: an example of
  {A}rtin-{M}umford.
\newblock In {\em The Brauer group (Sem., Les Plans-sur-Bex, 1980)}, volume 844
  of {\em Lecture Notes in Math.}, pages 103--128. Springer, Berlin, 1981.

\bibitem[Ume88]{Umemura-1988}
Hiroshi Umemura.
\newblock {Minimal rational threefolds. II.}
\newblock {\em Nagoya Math. J.}, 110:15--80, 1988.

\bibitem[Vos73]{Voskresenski-1973-re}
V.~E. Voskresenski{\v{i}}.
\newblock Fields of invariants of abelian groups.
\newblock {\em Uspehi Mat. Nauk}, 28:77--102, 1973.
\newblock Russian.

\bibitem[Web96]{Weber-1896}
H.~Weber.
\newblock {\em {Lehrbuch der Algebra. In zwei B�den. 2. Band.}}
\newblock {Braunschweig: F. Vieweg und Sohn. XIV + 796 S. $8^\circ$.}, 1896.

\end{thebibliography}
%\bibliographystyle{alpha}

\end{document}